\documentclass[11pt, a4paper]{article}
\usepackage{fullpage}

\usepackage[english]{babel}
\usepackage[utf8]{inputenc}

\usepackage{amsmath, amsfonts, amssymb}
\usepackage{mathtools}
\usepackage{bm}
\usepackage{comment}
\usepackage{authblk}
\usepackage{cite}
\usepackage[colorlinks=true, allcolors=blue]{hyperref}
\usepackage{graphicx}
\usepackage{subcaption}

\usepackage{booktabs}
\usepackage{makecell}

\usepackage[linesnumbered,ruled,vlined,plainruled]{algorithm2e}
\SetKwInput{KwInput}{input}

\SetCommentSty{mycommfont}

\DeclareMathOperator{\rank}{rank}

\DeclareMathOperator{\SVDr}{SVD_r}
\DeclareMathOperator{\svdr}{svdr}
\DeclareMathOperator{\QR}{QR}
\DeclareMathOperator{\HMT}{HMT}
\DeclareMathOperator{\Tropp}{Tropp}
\DeclareMathOperator{\HOSVD}{HOSVD}
\DeclareMathOperator{\STHOSVD}{STHOSVD}
\DeclareMathOperator{\TTSVD}{TTSVD}

\title{Low-rank nonnegative tensor approximation via alternating projections and sketching}
\author[1]{Azamat Sultonov}
\author[1,2]{Sergey Matveev}
\author[2]{Stanislav Budzinskiy}
\affil[1]{Faculty of Computational Mathematics and Cybernetics, Lomonosov Moscow State University, Moscow}
\affil[2]{Marchuk Institute of Numerical Mathematics RAS, Moscow}
\date{}

\begin{document}
\maketitle

\begin{abstract}
We show how to construct nonnegative low-rank approximations of nonnegative tensors in Tucker and tensor train formats. We use alternating projections between the nonnegative orthant and the set of low-rank tensors, using STHOSVD and TTSVD algorithms, respectively, and further accelerate the alternating projections using randomized sketching. The numerical experiments on both synthetic data and hyperspectral images show the decay of the negative elements and that the error of the resulting approximation is close to the initial error obtained with STHOSVD and TTSVD. The proposed method for the Tucker case is superior to the previous ones in terms of computational complexity and decay of negative elements. The tensor train case, to the best of our knowledge, has not been studied before.
\end{abstract}

\tableofcontents
\newpage

\section{Introduction}\label{sec:introduction}
Low-rank matrices and tensors have become ubiquitous in a multitude of tasks related to computational science and machine learning \cite{cichocki2017tensor, sidiropoulos2017tensor, khoromskij2018tensor}. These highly structured representations lead to significant reduction in storage and, simultaneously, to elegant and efficient numerical algorithms, which prove useful in solving high-dimensional ODEs and PDEs \cite{kazeev2014direct, matveev2016tensor, allmannrahn2022parallel}, recovering signals from scarce measurements \cite{liu2012tensor, budzinskiy2021tensor}, global optimization \cite{zheltkov2020global}, classification \cite{yang2017tensor}, etc.

In certain applications, the data are naturally nonnegative, and it is important to retain this property in the approximate representation. The data in question can be probability distributions \cite{dolgov2020approximation, novikov2021tensor, hur2022generative}, joint multidimensional concentrations of chemical compounds \cite{shcherbakova2019nonnegative1, manzini2021nonnegative}, multispectral images \cite{cichocki2009nonnegative}, audio \cite{leplat2020blind}, and others. The issue with the standard approaches to low-rank approximation, both for matrices and tensors, is that they give no guarantee of nonnegativity. 

An existing remedy is to employ the framework of nonnegative matrix factorization (NMF) \cite{gillis2020nonnegative} and nonnegative tensor factorization (NTF) \cite{cichocki2009nonnegative}. This is a collection of techniques, all of which are based on the idea that every individual factor of the low-rank decomposition should be nonnegative, thereby leading to a nonnegative approximation. The \textit{nonnegative} rank, however, can be larger than the usual rank of the matrix/tensor, so that the resulting representation becomes less compact (its further processing becomes less efficient too). 

An alternative view on the problem is less strict as it allows the factors of the low-rank decomposition to have negative entries \cite{vanluyten2008nonnegative}. The low-rank nonnegative matrix approximation (LRNMA) problem for $X \in \mathbb{R}^{m \times n}_{+}$ can be formulated as minimization of the Frobenius norm:
\begin{equation}
\label{eq:lrnma}
    \| X - Y \|_F \to \min \quad \text{s.t.} \quad Y \in \mathbb{R}^{m \times n}_{+}, \quad \rank\big( Y \big) \leq r.
\end{equation}
Some optimality properties of this best-approximation formulation were studied in \cite{grussler2015optimal}, viewed as a more general low-rank optimization problem with convex constraints \cite{grussler2018low, andersson2017convex}. Stepping away from the search for the \textit{best} approximation, a different approach based on alternating projections was proposed in \cite{song2020nonnegative} that provably locally converges to a \textit{good} approximation. Further developments in this direction concerned the computational efficiency of the alternating projections \cite{song2020tangent, matveev2022sketching} and an augmented Lagrangian method \cite{zhu2021approximate}. 

The most important difference between the algorithms for the NMF and LRNMA formulations is that in NMF, the low-rank iterates are guaranteed to be nonnegative; in LRNMA, however, the intermediate low-rank matrices \textit{do contain} negative elements, but they \textit{converge} to a low-rank nonnegative matrix as the number of negative elements and their magnitude decrease with iterations. This allows one to perceive the alternating projections approach to LRNMA as a filtering procedure that can reduce the magnitude of the negative entries of a low-rank matrix below a certain threshold, acceptable for a given application.

In this paper, we study how alternating projections can be used for multidimensional low-rank nonnegative tensor approximation (LRNTA) and accelerated with randomized low-rank projections, as was done for matrices in \cite{matveev2022sketching}. We focus on two popular tensor formats: the Tucker \cite{tucker1964extension} and tensor train (TT) \cite{oseledets2011tensor} decompositions. Tucker-LRNTA was the subject of \cite{jiang2020nonnegative}, where it was treated as consensus optimization between multiple LRNMA problems. To the best of our knowledge, TT-LRNTA is approached for the first time in our work.

\section{Preliminaries}
\subsection{Tensors}\label{sec:tensors}
We understand tensors as multidimensional arrays, i.e. elements of $\mathbb{R}^{n_1 \times \ldots \times n_d}$; see \cite{kolda2009tensor} for a gentle introduction. Matrices will be denoted by uppercase letters ($X, Y, \dots$) and tensors by bold uppercase letters ($\bm{X}, \bm{Y}, \dots$). The Frobenius norm of a tensor is the Euclidean norm induced by the standard inner product
\begin{equation*}
    \langle \bm{X}, \bm{Y} \rangle_F = \sum_{i_1, \ldots, i_d} \bm{X}(i_1, \ldots, i_d) \bm{Y}(i_1, \ldots, i_d), \quad \| \bm{X} \|_F = \sqrt{\langle \bm{X}, \bm{X} \rangle_F}, 
\end{equation*}
and the Chebyshev norm is defined as the entry of maximum magnitude
\begin{equation*}
    \| \bm{X} \|_C = \max_{i_1, \ldots, i_d} |\bm{X}(i_1, \ldots, i_d)|. 
\end{equation*}
The mode-$k$ unfolding of a $d$-dimensional tensor $\bm{X} \in \mathbb{R}^{n_1 \times \ldots \times n_d}$ for $1 \leq k \leq d$ is a matrix $X_{(k)} \in \mathbb{R}^{n_k \times \prod_{j \neq k} n_j}$ defined as
\begin{gather*}
    \bm{Y} = \textsc{permute}\Big(\bm{X}, \Big[k, 1, 2, \ldots, k-1, k+1, \ldots, d\Big]\Big), \\
    X_{(k)} = \textsc{reshape}\Big(\bm{Y}, \Big[n_k, \prod\nolimits_{j \neq k} n_j \Big] \Big),
\end{gather*}
with the help of Matlab-style $\textsc{permute}$ and $\textsc{reshape}$ operations. The $k$-th matricization of $\bm{X}$ for $1 \leq k \leq d-1$ is a matrix $X_{<k>} \in \mathbb{R}^{n_1 \ldots n_k \times n_{k+1} \ldots n_d}$ defined as
\begin{equation*}
    X_{<k>} = \textsc{reshape}\Big(\bm{X}, \Big[\prod\nolimits_{i = 1}^{k} n_i, \prod\nolimits_{j = k+1}^{d} n_j\Big]\Big).
\end{equation*}
The mode-$k$ product of a tensor $\bm{X} \in \mathbb{R}^{n_1 \times \ldots \times n_d}$ and a matrix $U \in \mathbb{R}^{m_k \times n_k}$ is a tensor 
\begin{equation*}
    \bm{Y} = \bm{X} \times_k U \in \mathbb{R}^{n_1 \times \ldots \times n_{k-1} \times m_k \times n_{k+1} \times \ldots \times n_d}, \quad Y_{(k)} = U X_{(k)}.
\end{equation*}

The Tucker, or multilinear, rank of a $d$-dimensional tensor $\bm{X}$ is an ordered $d$-tuple
\begin{equation*}
    \rank_{tucker}\big( \bm{X} \big) = \big( \rank\big( X_{(1)} \big), \ldots, \rank\big( X_{(d)} \big) \big).
\end{equation*}
If $\rank_{tucker}\big( \bm{X} \big) = (r_1, \ldots, r_d)$ then it is possible to represent $\bm{X}$ as a series of mode-$k$ products
\begin{equation*}
    \bm{X} = \bm{G} \times_1 U_1 \times_2 \ldots \times_d U_d
\end{equation*}
between a Tucker core $\bm{G} \in \mathbb{R}^{r_1 \times \ldots \times r_d}$ and Tucker factors $U_k \in \mathbb{R}^{n_k \times r_k}$. This representation is known as the Tucker decomposition of $\bm{X}$; it is not unique and none of the $r_k$ can be reduced without breaking the exact equality. For more details, see \cite{de2000multilinear, kolda2009tensor}.

The tensor train (TT) rank of a $d$-dimensional tensor $\bm{X}$ is an ordered $(d-1)$-tuple
\begin{equation*}
    \rank_{tt}\big( \bm{X} \big) = \big( \rank\big( X_{<1>} \big), \ldots, \rank\big( X_{<d-1>} \big) \big).
\end{equation*}
Let $\rank_{tt}\big( \bm{X} \big) = (r_1, \ldots, r_{d-1})$; then there exist two matrices $G_1 \in \mathbb{R}^{n_1 \times r_1}$, $G_d \in \mathbb{R}^{r_{d-1} \times n_d}$ and $d-2$ three-dimensional tensors $\bm{G}_k \in \mathbb{R}^{r_{k-1} \times n_k \times r_k}$ for $2 \leq k \leq d-2$ such that every entry of $\bm{X}$ can be calculated as
\begin{equation*}
    \bm{X}(i_1, \ldots, i_d) = \sum_{\alpha_1, \ldots, \alpha_{d-1}}G_1(i_1, \alpha_1) \bm{G}_2(\alpha_1, i_2, \alpha_2) \ldots \bm{G}_{d-1}(\alpha_{d-2}, i_{d-1}, \alpha_{d-1}) G_d(\alpha_{d-1}, i_d),
\end{equation*}
where each $\alpha_k$ ranges from $1$ to $r_k$. This is called the TT decomposition of $\bm{X}$ and $G_1$, $\{ \bm{G}_k \}$, $G_d$ are called TT cores. Like Tucker decomposition, it is not unique and exactness cannot be achieved with smaller values of $r_k$. Find more in \cite{oseledets2011tensor}.

\subsection{Sketching}\label{sec:sketching}
To compute the singular value decomposition (SVD) of an $m \times n$ matrix $X$ ($m \geq n$) is a classic problem in numerical linear algebra. The state-of-the-art algorithms (see \cite{dongarra2018singular} for a thorough review) require $O(m n^2)$ flops, whether you need all $n$ singular vectors or only $r \ll n$ dominant ones. This comes out as a serious bottleneck for computing low-rank approximations of large matrices; notably, matrices of size $n \times n^{d-1}$ with $d \geq 3$ appear in the HOSVD \cite{de2000multilinear} and TTSVD \cite{oseledets2009breaking} algorithms for tensor approximation in Tucker and TT formats, respectively. 

Randomized sketching is a powerful modern technique \cite{martinsson2020randomized}, which allows one to find \textit{good} (rather than \textit{best}) rank-$r$ approximations directly, using $O(mnk)$ flops with $k \geq r$, bypassing the full SVD. The main principle is to obtain a \textit{sketch} of the original matrix (by multiplying it with an $n \times k$ random matrix $\Psi$) that is smaller in size, yet its column space approximates the dominant singular subspace.

A better subspace can be obtained with the \textit{randomized subspace iteration} algorithm \cite[Alg. 4.4]{halko2011finding}, which computes the range of $(XX^T)^p X\Psi$ with $p \geq 0$. Practical implementation is presented in Alg.~\ref{rand_sub_iter}. The projection onto the corresponding orthonormal basis followed by truncated SVD (we denote it by $\SVDr$, where $r$ is the truncation rank) gives the randomized SVD algorithm \cite[Alg. 5.1]{halko2011finding} of the desired complexity (see Alg.~\ref{hmt}).

\begin{algorithm}[H]
\DontPrintSemicolon
  \KwInput{Data matrix $X \in \mathbb{R}^{m \times n}$, range sketch size $k \geq r$, number of iterations $p$, random matrix generator TestMatrix}
  $\Psi \leftarrow \text{TestMatrix}(n, k) \in \mathbb{R}^{n \times k}$\\
  $Z_1 \leftarrow X\Psi \in \mathbb{R}^{m \times k}$\\
  $[Q, R] \leftarrow \QR(Z_1)$\\
  \For {$j = 1, \dots, p$} {
      $Z_2 \leftarrow Q^TX \in \mathbb{R}^{k \times n}$\\
      $[Q, R] \leftarrow \QR(Z_2^T)$\\
      $Z_1 \leftarrow XQ \in \mathbb{R}^{m \times k}$\\
      $[Q, R] \leftarrow \QR(Z_1)$\\
  }
  \Return $Q$
\caption{Randomized subspace iteration \cite[Alg. 4.4]{halko2011finding}}
\label{rand_sub_iter}
\end{algorithm}

\begin{algorithm}[H]
\DontPrintSemicolon
  \KwInput{Data matrix $X \in \mathbb{R}^{m \times n}$, target rank $r$, estimated orthonormal basis $Q \in \mathbb{R}^{m \times k}$}
  $Z \leftarrow Q^T X \in \mathbb{R}^{k \times n}$\\
  $[U_r, \Sigma_r, V_r] \leftarrow \SVDr(Z)$\\
  $U_r \leftarrow Q U_r$\\
  \Return $U_{r}, \Sigma_{r}, V^T_{r}$
\caption{Randomized truncated SVD \cite[Alg. 5.1]{halko2011finding} (HMT)}
\label{hmt}
\end{algorithm}
Another approach to randomized SVD \cite{tropp2017practical} relies on two random matrices $\Psi \in \mathbb{R}^{n \times k}$ and $\Phi \in \mathbb{R}^{l \times m}$. It achieves the same asymptotic complexity as Alg.~\ref{hmt} but with a potentially smaller constant \cite{matveev2022sketching}. The first matrix $\Psi$ is used to obtain the orthonormal basis $Q$. Then, instead of projecting onto $Q$ directly, the second matrix $\Phi$ is employed through the Moore-Penrose pseudoinverse $(\Phi Q)^\dag \Phi X \in \mathbb{R}^{k \times n}$ (see Alg.~\ref{tropp}).

\begin{algorithm}[H]
\DontPrintSemicolon
  \KwInput{Data matrix $X \in \mathbb{R}^{m \times n}$, target rank $r$, range sketch size $k \geq r$, co-range sketch size $l \geq k$, random matrix generator TestMatrix}
  $\Psi \leftarrow \text{TestMatrix}(n, k) \in \mathbb{R}^{n \times k}$\\
  $\Phi \leftarrow \text{TestMatrix}(l, m) \in \mathbb{R}^{l \times m}$\\
  $Z \leftarrow X\Psi \in \mathbb{R}^{m \times k}$\\
  $[Q, R] \leftarrow \QR(Z)$\\
  $W \leftarrow \Phi Q \in \mathbb{R}^{l \times k}$\\
  $[P, T] \leftarrow \QR(W)$\\
  $G \leftarrow T^{-1}P^T\Phi X \in \mathbb{R}^{k \times n}$\\
  $[U_r, \Sigma_r, V_r] \leftarrow \SVDr(G)$\\
  $U_r \leftarrow QU_r$\\
  \Return $U_{r}, \Sigma_{r}, V^T_{r}$
\caption{Randomized truncated SVD \cite{tropp2017practical} (Tropp)}
\label{tropp}
\end{algorithm}

Before either of the algorithms is applied, the distribution of $\Psi$ (and $\Phi$) must be specified. In this paper, we will use matrices $\Psi \in \mathbb{R}^{n \times k}$ with iid Rademacher entries:
\begin{equation*}
\psi_{ij} \sim \text{Rad}, \quad \psi_{ij} = 
 \begin{cases}
   1,  &\text{with probability 1/2}\\
   -1, &\text{with probability 1/2}
 \end{cases}
\end{equation*}
For brevity, we will write $\HMT(p, k)$ for the combination of Algs.~\ref{rand_sub_iter}-\ref{hmt} and $\Tropp(k, l)$ for Alg.~\ref{tropp}.

\section{Low-rank nonnegative tensor approximation}\label{sec:methods}
\subsection{Matrix case}
An alternating-projections-based approach to solve the LRNMA problem was proposed in \cite{song2020nonnegative}. Consider two sets: the nonnegative orthant $\mathbb{R}^{m \times n}_{+}$ and the set of low-rank matrices
\begin{equation*}
    \mathcal{M}_{\leq r} = \{ X \in \mathbb{R}^{m \times n} : \rank \big( X \big) \leq r \}.
\end{equation*}
Both of them are closed, so for every matrix $X$ there are best (with respect to the Frobenius norm) nonnegative $\Pi_{\mathbb{R}^{m \times n}_{+}} \big( X \big) \in \mathbb{R}^{m \times n}_{+}$ and low-rank $\Pi_{\mathcal{M}_{\leq r}} \big( X \big) \subset \mathcal{M}_{\leq r}$ approximations, respectively. The former is unique due to convexity and is given by
\begin{equation*}
    \Pi_{\mathbb{R}^{m \times n}_{+}} \big( X \big) = \max\big( X, 0 \big).
\end{equation*}
The projection onto $\mathcal{M}_{\leq r}$, however, is not unique in general, but every matrix in $\Pi_{\mathcal{M}_{\leq r}} \big( X \big)$ is obtained as a truncated SVD of $X$ \cite{GolubVanLoanMatrix2013}. Hence, with a slight abuse of notation, we will write
\begin{equation*}
    \Pi_{\mathcal{M}_{\leq r}} \big( X \big) = \SVDr \big( X \big).
\end{equation*}

Given a nonnegative matrix $X \in \mathbb{R}^{m \times n}_{+}$, the algorithm from \cite{song2020nonnegative} then iterates between the two sets as follows:
\begin{equation*}
    X^{(0)} \gets X, \quad X^{(2k+1)} \gets \Pi_{\mathcal{M}_{\leq r}} \big( X^{(2k)} \big), \quad X^{(2k)} \gets \Pi_{\mathbb{R}^{m \times n}_{+}} \big( X^{(2k-1)} \big).
\end{equation*}
Clearly, one can also start with a low-rank matrix and change the order of the projections. Further papers considered more efficient alternating projections that use approximate rank truncation based on tangent spaces \cite{song2020tangent} and randomized sketching \cite{matveev2022sketching}.

As we noted in the Introduction and as is seen directly from the algorithm, the intermediate matrices are not simultaneously low-rank and nonnegative: the even iterates are only nonnegative and the odd iterates are only low-rank. However, both of the subsequences converge to a low-rank nonnegative matrix. It was proved in \cite{song2020nonnegative} that if the initial matrix $X$ is sufficiently close to the intersection $\mathcal{M}_{\leq r} \cap \mathbb{R}^{m \times n}_{+}$, the iterates converge to a quasioptimal solution of the LRNMA problem \eqref{eq:lrnma}. We aim to draw from these ideas to present multidimensional extensions of the alternating projections approach to Tucker-LRNTA and TT-LRNTA.

\subsection{Tucker case}
As we discussed in Subsec.~\ref{sec:tensors}, every tensor $\bm{X} \in \mathbb{R}^{n_1  \times \ldots \times n_d}$ admits an exact Tucker decomposition
\begin{equation*}
    \bm{X} = \bm{G} \times_1 U_1 \times_2 \ldots \times_d U_d, \quad \bm{G} \in \mathbb{R}^{r_1 \times \ldots \times r_d}, \quad U_k \in \mathbb{R}^{n_k \times r_k},
\end{equation*}
with $\rank_{tucker}\big( \bm{X} \big) = (r_1, \ldots, r_d)$. In practice, one seeks an \textit{approximate} Tucker decomposition 
\begin{equation*}
    \bm{X} \approx \bm{Y} = \bm{G} \times_1 U_1 \times_2 \ldots \times_d U_d
\end{equation*}
of given rank $\bm{r} = (r_1, \ldots, r_d)$ or such that the approximation is accurate up to a given threshold
\begin{equation*}
    \| \bm{X} - \bm{Y} \|_F < \varepsilon \| \bm{X} \|_F.
\end{equation*}

Both the exact and approximate Tucker decompositions can be constructed with the \textit{higher-order} SVD (HOSVD) algorithm \cite{de2000multilinear}: it computes the SVDs of the mode-$k$ unfoldings $X_{(k)}$, chooses the Tucker factors $U_k$ as the left singular factors, and computes the Tucker core by orthogonal projection
\begin{equation*}
    \bm{G} = \bm{X} \times_1 U_1^T \times_2 \ldots \times_d U_d^T.
\end{equation*}
Unlike the truncated SVD for matrices, HOSVD does not lead to the \textit{optimal} Tucker approximation; however, it is guaranteed to construct a quasioptimal one:
\begin{equation*}
    \| \bm{X} - \HOSVD_{\bm{r}} \big(\bm{X}\big) \|_F \leq \sqrt{d} \min_{\bm{Y}} \| \bm{X} - \bm{Y} \|_F, \quad \rank_{tucker}\big( \bm{Y} \big) \preceq \bm{r}.
\end{equation*}

The \textit{sequentially truncated} HOSVD (STHOSVD) is a more computationally efficient procedure with similar approximation properties \cite{vannieuwenhoven2012new}. We provide its pseudocode in Alg.~\ref{sthosvd}. Randomization can be employed to accelerate STHOSVD even further (cf. \cite{che2019randomized, ahmadi2021randomized}). By substituting $\HMT$ or $\Tropp$ in place of $\SVDr$, we get a family of algorithms with different rank-truncation strategies defined by $\svdr \in \mathcal{F} = \{ \SVDr, \HMT, \Tropp \}$.

\begin{algorithm}[H]
\DontPrintSemicolon
    \KwInput{Data tensor $\bm{X} \in \mathbb{R}^{n_1 \times \ldots \times n_d}$, target Tucker rank $\bm{r} = (r_1, \ldots, r_d)$}
    $\bm{G} \gets \bm{X}$\\
    \For {$k = 1, 2, \dots, d$} {
    $[U_{r_k}, \Sigma_{r_k}, V^T_{r_k}] \gets \SVDr(G_{(k)}, r_k)$\\
    $U_k \gets U_{r_k} \in \mathbb{R}^{n_k \times r_k}$ \\
    $G_{(k)} \gets \Sigma_{r_k} V^T_{r_k} \in \mathbb{R}^{r_k \times r_1\dots r_{k-1} n_{k+1} \dots n_d}$\\ %\Comment{Note: overwriting $G_{(k)}$ overwrites $G$}\\
    }
    \Return $\bm{G}, U_1, U_2, \dots, U_d$
\caption{Sequentially truncated higher-order SVD \cite{vannieuwenhoven2012new} (STHOSVD)}
\label{sthosvd}
\end{algorithm}

We will use STHOSVD and its randomized variants as approximate projections onto the closed set of tensors with low Tucker rank
\begin{equation*}
    \mathcal{M}_{\preceq \bm{r}}^{tucker} = \{ \bm{X} \in \mathbb{R}^{n_1 \times \ldots \times n_d} : \rank_{tucker} \big( \bm{X} \big) \preceq \bm{r} \}.
\end{equation*}
Since the projection onto the nonnegative orthant $\mathbb{R}^{n_1 \times \ldots \times n_d}_{+}$ is the same as for matrices,
\begin{equation*}
    \Pi_{\mathbb{R}^{n_1 \times \ldots \times n_d}_{+}} \big( \bm{X} \big) = \max\big( \bm{X}, 0 \big),
\end{equation*}
we can formulate an alternating projection algorithm NSTHOSVD for the LRNTA problem in Tucker format; see Alg.~\ref{nsthosvd}.

\begin{algorithm}[H]
\DontPrintSemicolon
  \KwInput{Data tensor $\bm{X} \in \mathbb{R}^{n_1 \times \ldots \times n_d}$, target Tucker rank $\bm{r} = (r_1, \ldots, r_d)$, number of iterations $s$, rank-truncation strategy $\svdr \in \mathcal{F}$}
  $\bm{X}^{(0)} \gets \bm{X}$ \\
  \For {$i = 1, 2, \dots, s$}
    {
        $\bm{X}^{(i)} \gets \max(\bm{X}^{(i-1)},0)$ \\
        $[\bm{G}, U_1, \dots, U_d] \gets \STHOSVD(\bm{X}^{(i)}, \bm{r}, \svdr)$\\
        $\bm{X}^{(i)} \gets \bm{G} \times_1 U_1 \times_2 \ldots \times_d U_d \in \mathbb{R}^{n_1 \times \dots \times n_d}$\\
    }
  \Return $\bm{G}, U_1, U_2, \ldots, U_d$
\caption{STHOSVD-based alternating projections (NSTHOSVD)}
\label{nsthosvd}
\end{algorithm}

Note that a different alternating projections algorithm NLRT was proposed for Tucker-LRNTA in \cite{jiang2020nonnegative} and its convergence was proved. While we suggest to use approximate projections onto the actual sets of interest $\mathbb{R}^{n_1 \times \ldots \times n_d}_{+}$ and $\mathcal{M}_{\preceq \bm{r}}^{tucker}$, NLRT performs exact projections onto modified sets
\begin{gather*}
    \Omega_1 = \big\{ (\bm{X}_1, \ldots, \bm{X}_d) : \bm{X}_1 = \ldots = \bm{X}_d \in \mathbb{R}^{n_1 \times \ldots \times n_d}_{+} \big\}, \\
    \Omega_2 = \big\{ (\bm{X}_1, \ldots, \bm{X}_d) : \big(\bm{X}_k\big)_{(k)} \in \mathcal{M}_{\leq r_k},~1 \leq k \leq d  \big\},
\end{gather*}
which are subsets of the Cartesian product $\big( \mathbb{R}^{n_1 \times \ldots \times n_d} \big) \times \ldots \times \big( \mathbb{R}^{n_1 \times \ldots \times n_d} \big)$. This is essentially a consensus optimization problem, where LRNMA is computed for every mode-$k$ unfolding individually.

\subsection{Tensor train case}
In the same vein, exact and approximate TT decompositions of a tensor $\bm{X}$ can be computed with TTSVD \cite{oseledets2011tensor}, whose pseudocode we show in Alg.~\ref{ttsvd}. The resulting approximation of given TT rank $\bm{r} = (r_1, \ldots, r_{d-1})$ is quasioptimal
\begin{equation*}
    \| \bm{X} - \TTSVD_{\bm{r}} \big(\bm{X}\big) \|_F \leq \sqrt{d-1} \min_{\bm{Y}} \| \bm{X} - \bm{Y} \|_F, \quad \rank_{tt}\big( \bm{Y} \big) \preceq \bm{r},
\end{equation*}
and is typically used as an approximate projection onto the set of tensors with low TT rank
\begin{equation*}
    \mathcal{M}_{\preceq \bm{r}}^{tt} = \{ \bm{X} \in \mathbb{R}^{n_1 \times \ldots \times n_d} : \rank_{tt} \big( \bm{X} \big) \preceq \bm{r} \}.
\end{equation*}

\begin{algorithm}[H]
\DontPrintSemicolon
  \KwInput{Data tensor $\bm{X} \in \mathbb{R}^{n_1 \times \ldots \times n_d}$, target TT rank $\bm{r} = (r_1, \ldots, r_{d-1})$}
    $G \gets \textsc{reshape}\big(\bm{X}, [n_1, n_2 \ldots n_d]\big)] \in \mathbb{R}^{n_1 \times n_2 \ldots n_d}$ \\
    $[U_{r_1}, \Sigma_{r_1}, V^T_{r_1}] \gets \SVDr(G,\ r_1)$ \\
    $G_1 \gets U_{r_1}$ \\
    \For {$k = 2, \dots, d-1$} {
        $G \gets \textsc{reshape}(\Sigma_{r_{k-1}} V^T_{r_{k-1}},\ [r_{k-1} n_{k}, n_{k+1} \ldots n_d]) \in \mathbb{R}^{r_{k-1} n_{k} \times n_{k+1} \ldots n_{d}}$ \\
        $[U_{r_k}, \Sigma_{r_k}, V^T_{r_k}] \gets \SVDr(G,\ r_k)$ \\
        $\bm{G}_k \gets \textsc{reshape}\big(U_{r_k}, [r_{k-1}, n_k, r_k]) \in \mathbb{R}^{r_{k-1} \times n_k \times r_k}$ \\
    }
    $G_d \gets \Sigma_{r_{d-1}} V^T_{r_{d-1}}$ \\
    \Return $G_1, \bm{G}_2, \ldots, \bm{G}_{d-1}, G_{d}$
\caption{TTSVD \cite{oseledets2011tensor}}
\label{ttsvd}
\end{algorithm}

Swapping $\SVDr$ for a randomized SVD algorithm such as $\HMT$ or $\Tropp$ is a straightforward way to reduce the computational complexity (cf. \cite{huber2017randomized, che2019randomized, kressner2022streaming}). We use these variants of TTSVD to solve the TT-LRNTA problem with alternating projections. Our algorithm NTTSVD is presented in Alg.~\ref{nttsvd}.

\begin{algorithm}[H]
\DontPrintSemicolon
  \KwInput{Data tensor $\bm{X} \in \mathbb{R}^{n_1 \times \ldots \times n_d}$, target TT rank $\bm{r} = (r_1, \ldots, r_{d-1})$, number of iterations $s$, rank-truncation strategy $\svdr \in \mathcal{F}$}
  $\bm{X}^{(0)} \gets \bm{X}$ \\
  \For {$i = 1, 2, \dots, s$}
    {
        $\bm{X}^{(i)} \gets \max(\bm{X}^{(i-1)},0)$ \\
        $[G_1, \bm{G}_2, \ldots, \bm{G}_{d-1}, G_d] \gets \TTSVD(\bm{X}^{(i)}, \bm{r}, \svdr)$\\
        $\bm{X}^{(i)} \gets \sum\limits_{\alpha_1, \ldots, \alpha_{d-1}}G_1(j_1, \alpha_1) \bm{G}_2(\alpha_1, j_2, \alpha_2) \ldots \bm{G}_{d-1}(\alpha_{d-2}, j_{d-1}, \alpha_{d-1}) G_d(\alpha_{d-1}, j_d)$\\
    }
    \Return $G_1, \bm{G}_2, \ldots, \bm{G}_{d-1}, G_d$
\caption{TTSVD-based alternating projections (NTTSVD)}
\label{nttsvd}
\end{algorithm}

\subsection{Computational complexities}\label{subsec:complexity}
In this section, we will estimate the one-iteration complexity of NSTHOSVD and NTTSVD. The most computationally expensive part of both algorithms is rank truncation, for which we have 3 options: SVD, $\HMT(p, k)$, and $\Tropp(k, l)$. For an $m \times n$ matrix with $m \geq n$ and truncation rank $r$, the standard truncated SVD costs $O(mn^2)$ floating point operations \cite{GolubVanLoanMatrix2013, dongarra2018singular}. As for the randomized approaches, $\HMT(p, k)$ requires $O(mn (pk+r))$ operations and $\Tropp(k, l)$ runs in $O(mn(r+k+l))$. While the asymptotic complexities are independent of the random matrix distribution, the constants can be reduced if structured random matrices are used (see \cite{matveev2022sketching}).

Turning to STHOSVD and TTSVD, let us denote by $r$ the maximum truncation rank of the tensor (Tucker or TT), let $n = \max(n_1, \ldots, n_d)$, and assume that $r$ is small compared to $n$. Every step of STHOSVD (Alg.~\ref{sthosvd}) is dominated by rank truncation of unfolding matrices, whose sizes decrease progressively from $n \times n^{d-1}$ down to $n \times r^{d-1}$. The whole procedure then costs $O(n^{d+1})$ for SVD, $O(n^d (pk+r))$ for $\HMT(p, k)$, and $O(n^d (r + k + l))$ for $\Tropp(k, l)$. The analysis of TTSVD (Alg.~\ref{ttsvd}) is absolutely analogous: the matricizations become smaller as $n \times n^{d-1}$, $rn \times n^{d-2}$, \dots, $rn \times n$, so that the total cost is also $O(n^{d+1})$, $O(n^d (pk+r))$, and $O(n^d (r + k + l))$, respectively.

To compute a single element of a tensor given in the Tucker and TT formats, $O(d r^2)$ and $O(r^d)$ operations are needed, respectively. However, given the rich multilinear structure, it is possible to build the full tensor directly, using $O(n^d r)$ in both cases, which is faster than computing $n^d$ elements individually. The nonnegative projections then require $O(n^d r)$ operations. This tells us that the overall complexity of a single iteration of NSTHOSVD and NTTSVD is defined by rank truncation, and we list it in Tab.~\ref{complexity}. With randomization, we reach balance in complexity of low-rank and nonnegative projections: both of them scale linearly with $n^d$, the number of elements of the tensor.
\begin{table}[h]\centering
\begin{tabular}{@{}lccc@{}}
    \toprule
     & SVD & $\HMT(p, k)$ & $\Tropp(k, l)$ \\
    \midrule
    NSTHOSVD & $O(n^{d+1})$ & $O(n^d (pk + r))$ & $O(n^d (r + k + l))$ \\
    NTTSVD & $O(n^{d+1})$ & $O(n^d (pk + r))$ & $O(n^d (r + k + l))$ \\
    \bottomrule\\
\end{tabular}
\caption{The asymptotic computational complexity of a single iteration of NSTHOSVD and NTTSVD algorithms with different rank-truncation methods.}
\label{complexity}
\end{table}

\subsection{Related work: existing theoretical analysis of alternating projections}
\indent Originally, the method of alternating projections was developed to compute the orthogonal projection of a given point onto the intersection of two closed subspaces in a Hilbert space \cite{escalante2011alternating}. Provided that the sum of these subpsaces is closed too, the iterates converge linearly in norm for any starting point with a rate defined by the Friedrichs angle between the two subspaces \cite{deutsch1984rate, kayalar1988error}.\\
\indent One way to generalize closed subspaces is to consider closed convex sets (such as the nonnegative orthant $\mathbb{R}^{n_1 \times \ldots \times n_d}_{+}$). For a pair of two closed convex sets, the method of alternating projections no longer finds the best approximation of the starting point, but solves the feasibility problem: converges to an arbitrary point in the intersection \cite{bauschke1994dykstra}. The iterates converge for every starting point, but do so only weakly. To prove strong convergence or even linear convergence, it is required that the pair of closed convex sets satisfies certain regularity assumptions \cite{bauschke1993convergence}.\\
\indent The setting of closed convex sets is the most natural for the method of alternating projections. Indeed, owing to the Bunt-Motzkin theorem, every point of a finite-dimensional Euclidean space has a unique best approximation by a set if and only if this set is closed and convex \cite{deutsch2001best}. Meanwhile, alternating projections have been extensively used for nonconvex sets as well, even though the projections are not uniquely defined for them (which complicates the global convergence analysis). Many sets that appear in practice are, however, \textit{prox-regular}: the projections onto them are locally unique, i.e. the best approximation is uniquely defined for every point that is close enough to the set \cite{poliquin2000local}. This property is shared, for example, by closed convex sets (obviously) and smooth manifolds \cite{lewis2008alternating}. When at least one of the two sets is prox-regular and their intersection satisfies certain regularity properties, the method of alternating projections locally linearly converges to a point in the intersection \cite{lewis2009local}. Notably, if the two sets are smooth manifolds, the regularity assumption translates to their intersection being transversal \cite{lewis2008alternating} and can be further relaxed to nontangential intersections \cite{andersson2013alternating}. Moreover, it was shown in \cite{andersson2013alternating} that alternating projections on smooth manifolds converge to quasioptimal approximations.\\
\indent The sets of rank-$r$ matrices, tensors of Tucker rank $\bm{r}$, and tensors of TT rank $\bm{r}$ are smooth manifolds \cite{uschmajew2020geometric} and the set of low-rank matrices $\mathcal{M}_{\leq r}$ is prox-regular at every rank-$r$ matrix \cite{luke2013prox}. This suggests good behavior of the alternating projections for the LRNMA and LRNTA problems. The matrix case, where truncated SVD produces the optimal low-rank approximation, was theoretically addressed in \cite{song2020nonnegative}. In the tensor case, STHOSVD and TTSVD lead to quasioptimal projections, which will make the convergence analysis for NSTHOSVD and NTTSVD more delicate.\\
\indent In the present work, we do not concentrate on the theoretical side of why NSTHOSVD and NTTSVD converge to low-rank nonnegative tensors. Instead, we rely on strong numerical evidence showing that they do in a number of different experiments. The successful outcome, in turn, motivates us to prove rigorous convergence guarantees in the future papers.

\section{Numerical experiments}\label{sec:numerical}
In this section, we evaluate and compare the performance of deterministic and randomized variants of the NSTHOSVD and the NTTSVD algorithms, which were introduced in Sec.~\ref{sec:methods}. The examples we consider are 
\begin{itemize}
    \item the Hilbert tensor (Subsec.~\ref{subsec:hilbert}),
    \item a mixture of multidimensional Gaussians (Subsec.~\ref{subsec:gaussian}),
    \item a hyperspectral image (Subsec.~\ref{subsec:hyperspectral}).
\end{itemize}
For every Tucker-LRNTA experiment, we provide the respective results obtained with the NLRT algorithm \cite{jiang2020nonnegative} to compare its performance with NSTHOSVD. Note that comparing NTTSVD with NSTHOSVD and NLRT is not particularly meaningful since the difference in running times are mostly dictated by how well a given dataset is approximated in Tucker and TT tensor formats.

All the experiments were carried out in Python (3.9.12) with Intel(R) Core(TM) i3-8130U CPU@2.20GHz and 8GB of RAM.`

\subsection{Hilbert tensor}\label{subsec:hilbert}
Our first example of an approximately low-rank nonnegative tensor is the Hilbert tensor
\begin{equation*}
    \bm{X}(i_1, \ldots, i_d) = \frac{1}{i_1+\ldots+i_d - d  + 1},
\end{equation*}
which is a multidimensional extension of the well-known Hilbert matrix. In Table~\ref{tab:Hilbert_timinigs_errors}, we present the results obtained with NSTHOSVD, NTTSVD, and NLRT when applied to a 3-dimensional Hilbert tensor of size $128 \times 128 \times 128$ with Tucker ranks $\bm{r} = (3, 2, 4)$ and TT-ranks $\bm{r} = (3, 2)$.

\begin{figure}[t]
\centering
\begin{subfigure}[b]{0.45\textwidth}
	\includegraphics[width=1.0\textwidth]{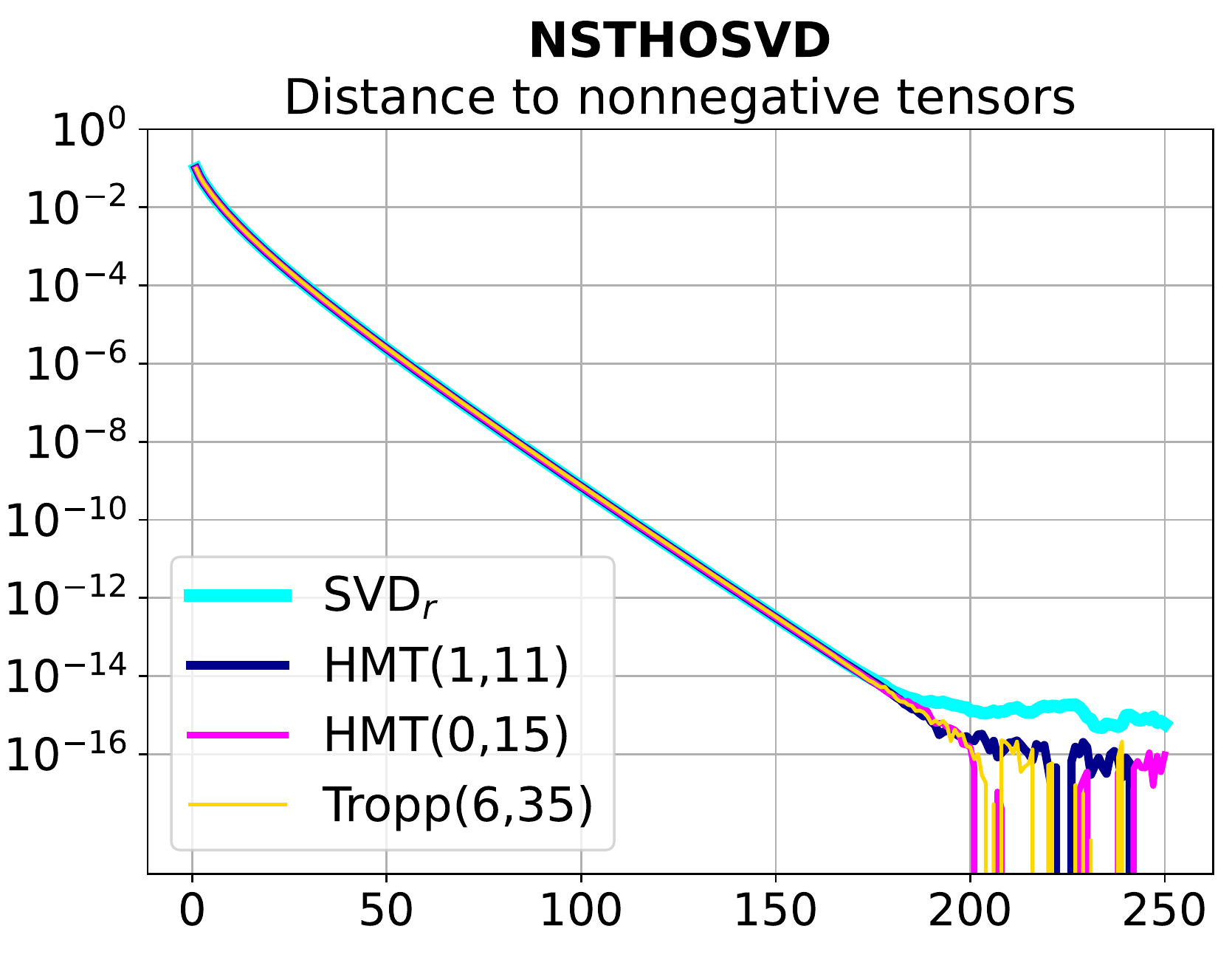}
	\caption{}
\end{subfigure}\hfill%
\begin{subfigure}[b]{0.45\textwidth}
	\includegraphics[width=1.0\textwidth]{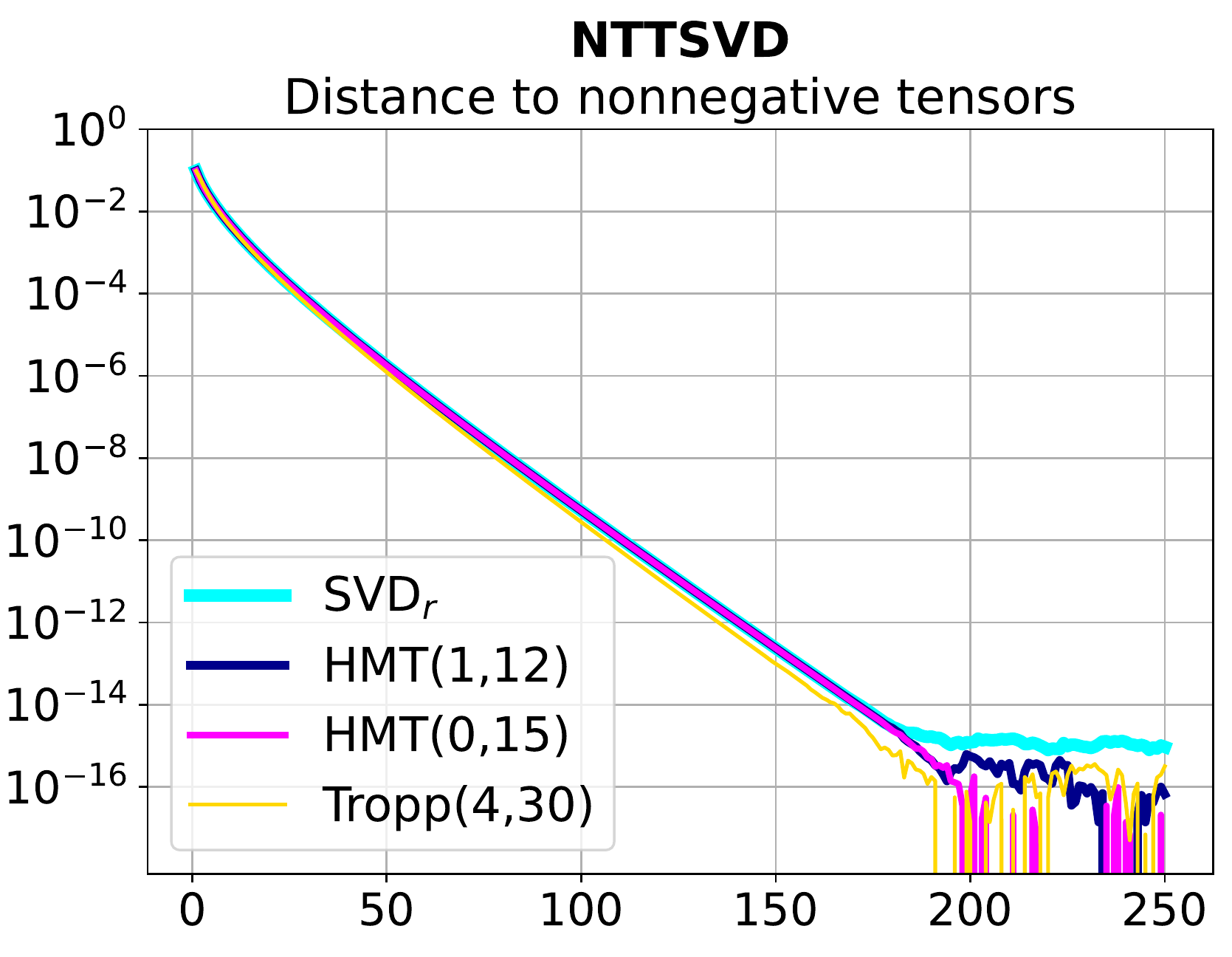}
	\caption{}
\end{subfigure}
\caption{Comparison of deterministic and randomized variants of NSTHOSVD~(a) and NTTSVD~(b) for low-rank nonnegative tensor approximation of a $128 \times 128 \times 128$ Hilbert tensor with Tucker ranks $(3, 2, 4)$ and TT-ranks $(3, 2)$: the Frobenius norm of the negative part over 250 iterations.}
\label{fig:Hilbert_convergence}
\end{figure}

First of all, note that the low-rank approximations obtained with simple STHOSVD and TTSVD contain negative elements whose total Frobenius norm is about $10^{-2}$. Using 250 iterations of NSTHOSVD and NTTSVD (deterministic or randomized), we can reduce their Frobenius norm down to 5 double-precision machine epsilons. The total number of negative elements also decreases, and sometimes we manage to remove them completely. Remarkably, the relative approximation errors grow by only about 3\% (Frobenius) and 7.5\% (Chebyshev) so that the resulting low-rank tensors are still \textit{good} approximations. The running time is an important aspect too: using randomized sketching, we achieved speed-up factors of about 7-11, compared to NSTHOSVD/NTTSVD with deterministic low-rank projections. All the variants show identical linear decay of the Frobenius norm of the negative elements (see Fig.~\ref{fig:Hilbert_convergence}).

Before comparing NSTHOSVD with NLRT, we would like to point out that NLRT never actually forms a tensor in the Tucker format. At every iteration, it computes low-rank approximations $(\bm{X}_1, \ldots, \bm{X}_d) \in \Omega_2$ of every unfolding (like HOSVD) but proceeds to the nonnegative projection without computing the Tucker core (unlike HOSVD). While NLRT iterations (all its $d$ components) do converge to a single nonnegative tensor with low Tucker ranks, it is not formed explicitly. In Table~\ref{tab:Hilbert_timinigs_errors}, the NLRT-related results on relative errors and negative elements are computed based on the following low-rank STHOSVD approximation,
\begin{equation}
\label{eq:nlrt_aux}
    \hat{\bm{X}} = \STHOSVD_{\bm{r}}\Big( \frac{1}{d} \sum\nolimits_{k = 1}^{d} \Pi_{\mathbb{R}^{n_1 \times \ldots \times n_d}_{+}} \big( \bm{X}_k \big) \Big),
\end{equation}
and the running time is shown for the original NLRT (i.e. excluding the additional STHOSVD).

The results in Table~\ref{tab:Hilbert_timinigs_errors} show that NSTHOSVD is superior to NLRT in terms of speed and mitigating negative elements, and they achieve similar relative errors. Deterministic NSTHOSVD is about 3 times faster than NLRT (just as STHOSVD is asymptotically $d$ times faster than HOSVD) and its randomized versions are 23-39 times faster. The Frobenius norm of the negative elements of $\hat{\bm{X}}$ is of order $10^{-9}$, which is an improvement over the initial $10^{-2}$, but is much larger than what NSTHOSVD achieves. We also compare the properties of $\bm{X}_1$, $\bm{X}_2$, $\bm{X}_3$, and $\hat{\bm{X}}$ after 250 iterations of NLRT (see Table~\ref{tab:Hilbert_NLRT_details}). 

\begin{table}[h]\centering
\begin{tabular}{@{}lccccc@{}}
    \toprule
    Method & \thead{Running time\\(Second)} & \thead{Relative error\\(Frobenius)} & \thead{Relative error\\(Chebyshev)} & \thead{Negative\\elements\\(Frobenius)} & \thead{Negative\\elements\\(\%)}\\ \hline
    TTSVD                   & 0.6 & $7.72\cdot 10^{-2}$ & $3.67\cdot 10^{-1}$ & $9.7\cdot 10^{-3}$ & $6.3\cdot 10^{-3}$ \\ \hline
    NTTSVD, SVD$_r$         & 121 & $7.88\cdot 10^{-2}$ & $3.94\cdot 10^{-1}$ & $9.3\cdot 10^{-16}$ & $1.9\cdot 10^{-4}$\\
    NTTSVD, HMT(1, 12)      & 16  & $7.88\cdot 10^{-2}$ & $3.94\cdot 10^{-1}$ & 0.0 & 0.0\\
    NTTSVD, HMT(0, 15)      & 13  & $7.88\cdot 10^{-2}$ & $3.94\cdot 10^{-1}$ & $1.2\cdot 10^{-16}$ & $9.5\cdot 10^{-5}$\\
    NTTSVD, Tropp(4, 30)    & 10  & $8.13\cdot 10^{-2}$ & $3.82\cdot 10^{-1}$ & 0.0 & 0.0\\\hline
    STHOSVD                 & 0.5 & $7.72\cdot 10^{-2}$ & $3.67\cdot 10^{-1}$ & $9.7\cdot 10^{-3}$ & $6.3\cdot 10^{-3}$ \\ \hline
    NSTHOSVD, SVD$_r$       & 118 & $7.89\cdot 10^{-2}$ & $3.95\cdot 10^{-1}$ & $5.9\cdot 10^{-16}$ & $1.4\cdot 10^{-4}$\\
    NSTHOSVD, HMT(1, 11)    & 17  & $7.89\cdot 10^{-2}$ & $3.95\cdot 10^{-1}$ & 0.0 & 0.0\\
    NSTHOSVD, HMT(0, 15)    & 14  & $7.89\cdot 10^{-2}$ & $3.95\cdot 10^{-1}$ & 0.0 & 0.0\\
    NSTHOSVD, Tropp(6, 35)  & 10  & $7.89\cdot 10^{-2}$ & $3.94\cdot 10^{-1}$ & $2.5\cdot 10^{-16}$ & $4.8\cdot 10^{-5}$\\ \hline
    NLRT          & 390 & $7.88\cdot 10^{-2}$ & $3.99\cdot 10^{-1}$ & $8.6\cdot 10^{-10}$ & $3.3\cdot 10^{-4}$\\
    \bottomrule
\end{tabular}
\caption{Comparison of NTTSVD, NSTHOSVD, and NLRT for low-rank nonnegative tensor approximation of a $128 \times 128 \times 128$ Hilbert tensor with Tucker ranks $(3, 2, 4)$ and TT-ranks $(3, 2)$: running times, relative errors, and negative elements after 250 iterations.}
\label{tab:Hilbert_timinigs_errors}
\end{table}

\begin{table}[h]\centering
\begin{tabular}{@{}cccccc@{}}
    \toprule
    Tensor & \thead{Relative\\error\\(Frobenius)} & \thead{Relative\\error\\(Chebyshev)} & \thead{Negative\\elements\\(Frobenius)} & \thead{Negative\\elements\\(Chebyshev)} & \thead{Negative\\elements\\(\%)}\\ \hline
    $\bm{X}_1$                  & $7.88\cdot 10^{-2}$ & $3.99\cdot 10^{-1}$ & $1.9\cdot 10^{-10}$ & $1.1\cdot 10^{-10}$ & $3.3\cdot 10^{-4}$ \\
    $\bm{X}_2$                  & $7.88\cdot 10^{-2}$ & $3.99\cdot 10^{-1}$ & $8.5\cdot 10^{-10}$ & $5.1\cdot 10^{-10}$ & $3.3\cdot 10^{-4}$ \\
    $\bm{X}_3$                  & $7.88\cdot 10^{-2}$ & $3.99\cdot 10^{-1}$ & $4.6\cdot 10^{-11}$ & $2.7\cdot 10^{-11}$ & $3.3\cdot 10^{-4}$ \\
    $\hat{\bm{X}}$ & $7.88\cdot 10^{-2}$ & $3.99\cdot 10^{-1}$ & $8.6\cdot 10^{-10}$ & $5.2\cdot 10^{-10}$ & $3.3\cdot 10^{-4}$ \\
    \bottomrule
\end{tabular}
\caption{Comparison of the NLRT components $\{ \bm{X}_k \}_{k = 1}^{3}$ and the auxiliary tensor $\hat{\bm{X}}$ in Tucker format for low-rank nonnegative tensor approximation of a $128 \times 128 \times 128$ Hilbert tensor with Tucker ranks $(3, 2, 4)$: relative errors and negative elements after 250 iterations.}
\label{tab:Hilbert_NLRT_details}
\end{table}

\subsection{Multidimensional Gaussian mixture}\label{subsec:gaussian}
In the second experiment, we test our approach on synthetic data that are an example of the low-rank density approximation problem: a multidimensional mixture of Gaussians
\begin{equation*}
    f(x) = \sum_{j = 1}^{m} \alpha_j \exp\Big( (x - \mu_j)^\top A_j^{-1} (x - \mu_j) \Big), \quad x \in \mathbb{R}^{d},
\end{equation*}
with weights $\alpha_j \in \mathbb{R}$, means $\mu_j \in \mathbb{R}^{d}$, and covariance matrices $A_j \in \mathbb{R}^{d \times d}$. Every individual Gaussian has approximately low rank (it is a rank-1 tensor if $A_j$ is diagonal) so the mixture can be approximated as well. We consider the mixture $f(x)$ in a hypercube $[-a, a]^d$ and discretize the domain on an equidistant tensor-product grid with step $2a / (n - 1)$, which gives a $d$-dimensional $n \times \ldots \times n$ tensor $\bm{X}$.

We choose the 4-dimensional scenario with a balanced mixture of 2 Gaussians ($\alpha_1 = \alpha_2$) with the following means,
\begin{equation*}
    \mu_1 = 
    \begin{bmatrix}
        0 & 0 & 0 & 0
    \end{bmatrix}^\top, \quad
    \mu_1 = 
    \begin{bmatrix}
        0.5 & -0.5 & 0.5 & -0.5
    \end{bmatrix}^\top,
\end{equation*}
and covariance matrices,
\begin{equation*}
    A_1 = 
    \begin{bmatrix}
        0.403 & 0.236 & 0.159 & 0.188 \\
        0.236 & 0.422 & 0.193 & 0.313 \\
        0.159 & 0.193 & 0.124 & 0.164 \\
        0.188 & 0.313 & 0.164 & 0.288
    \end{bmatrix}, \quad
    A_2 = 
    \begin{bmatrix}
        0.173 & 0.229 & 0.200 & 0.191 \\
        0.229 & 0.347 & 0.254 & 0.201 \\
        0.200 & 0.254 & 0.348 & 0.252 \\
        0.191 & 0.201 & 0.252 & 0.360
    \end{bmatrix}. 
\end{equation*}
We take $a = 1$ as the size of the domain and $n = 64$ as the size of the tensor.

\begin{figure}[h]
\begin{subfigure}[b]{1.0\linewidth}
\centering
	\includegraphics[width=0.45\linewidth]{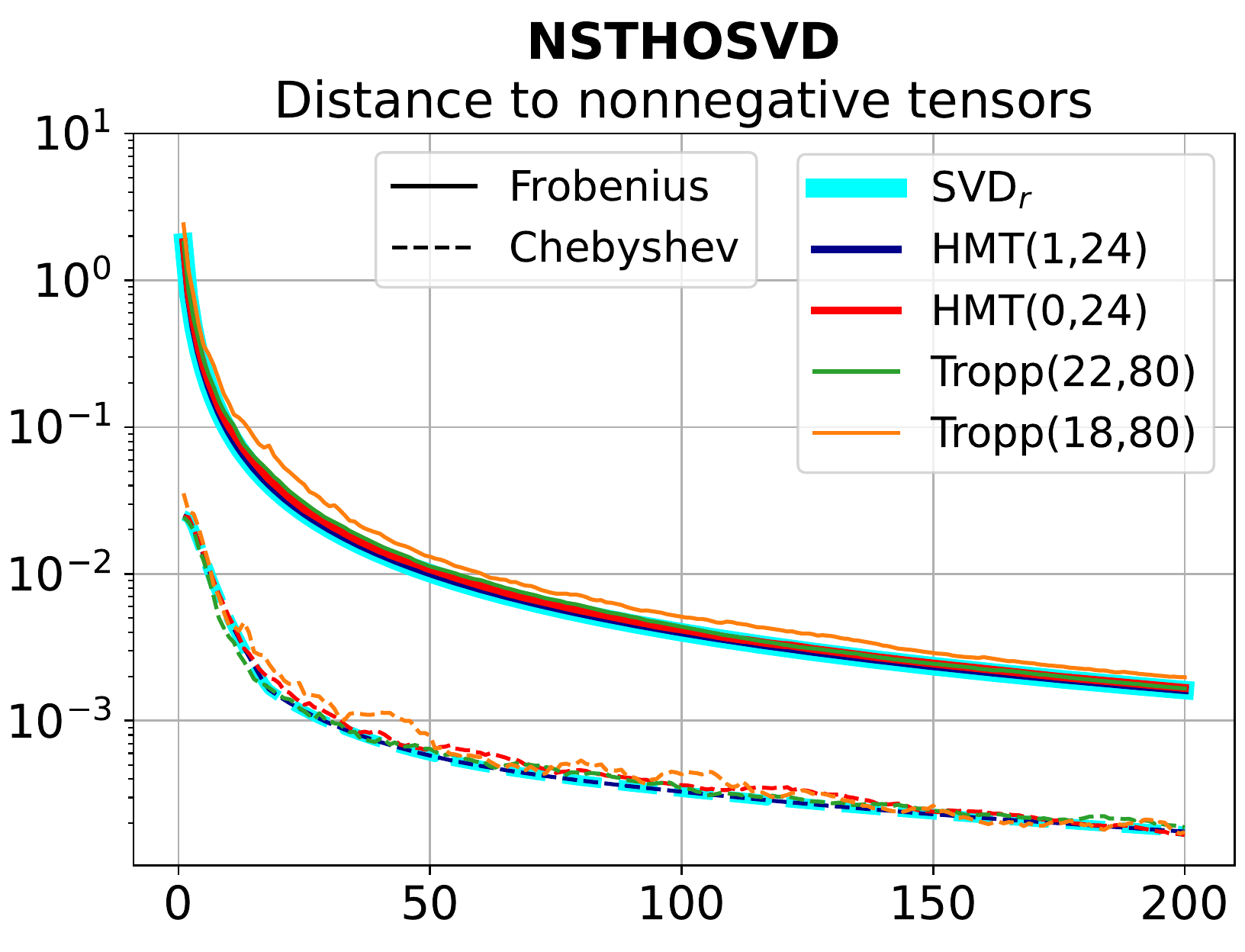} \hfill%
	\includegraphics[width=0.45\linewidth]{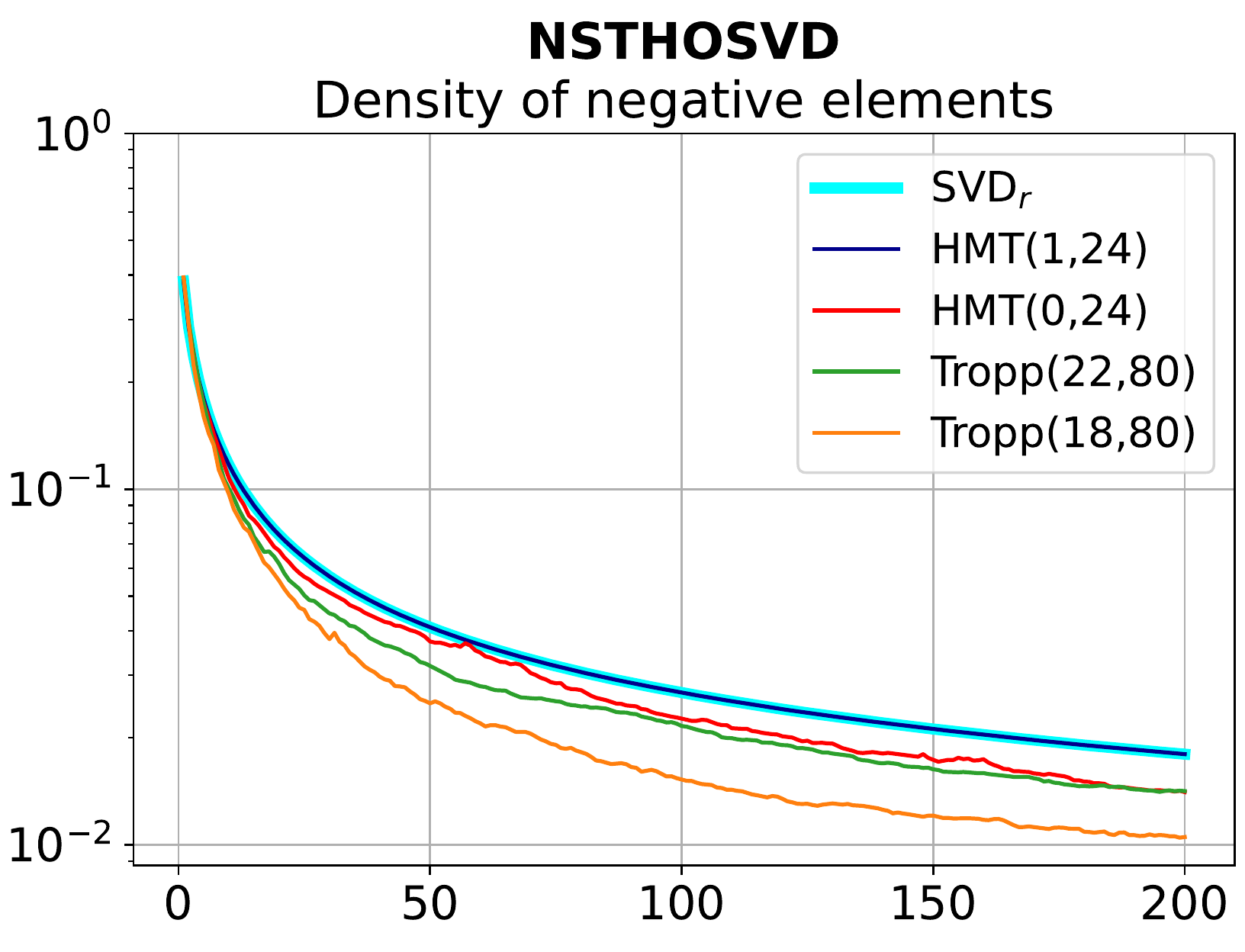}
	\caption{}
\end{subfigure}
\begin{subfigure}[b]{1.0\linewidth}
\centering
	\includegraphics[width=0.45\linewidth]{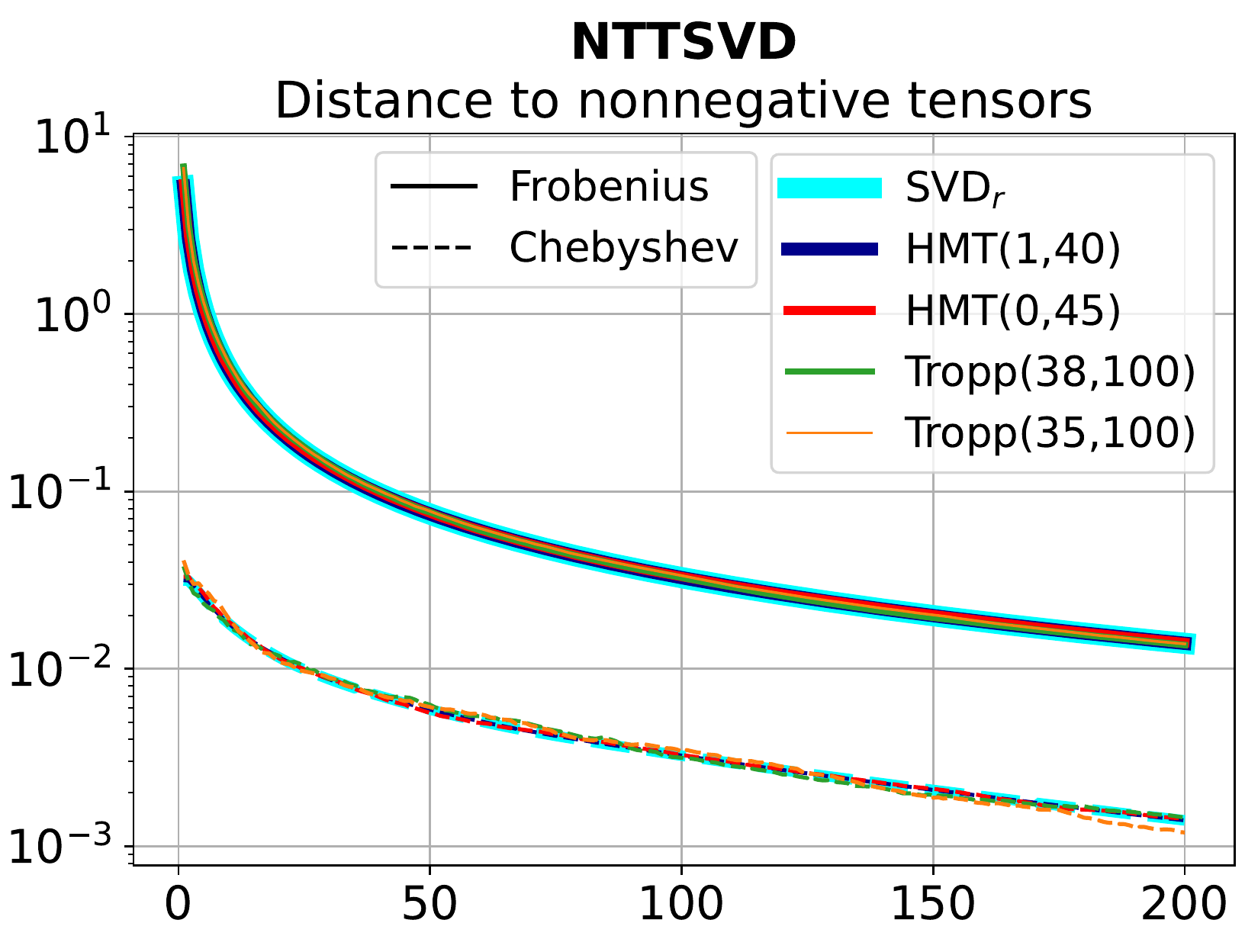} \hfill%
	\includegraphics[width=0.45\linewidth]{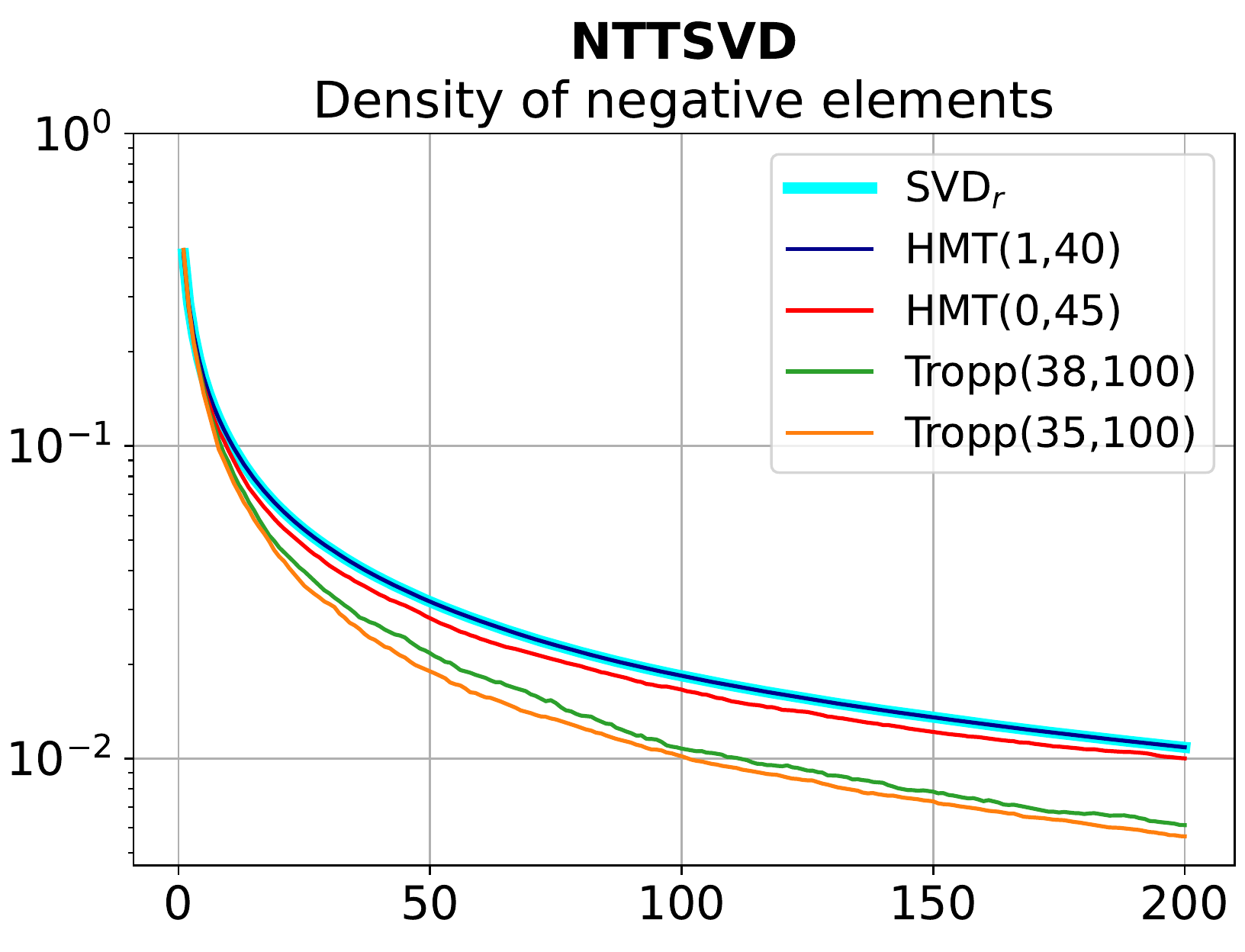}
	\caption{}
\end{subfigure}
\caption{Comparison of deterministic and randomized variants of NSTHOSVD~(a) and NTTSVD~(b) for low-rank nonnegative tensor approximation of a $64 \times 64 \times 64 \times 64$ Gaussian mixture with Tucker ranks $(14,14,14,14)$ and TT-ranks $(10,20,10)$: the Frobenius and Chebyshev norms of the negative part and the density of negative elements over 200 iterations.}
\label{fig:Gaussiance_convergence_rate}
\end{figure}

We carried out 200 iterations of NTTSVD, NSTHOSVD, and NLRT with Tucker ranks $(14,14,14,14)$ and TT-ranks $(10,20,10)$; see Table~\ref{tab:Gaussian_timings_errors}. Simple STHOSVD and TTSVD produce tensors with many, about 40\%, negative elements. The application of NTTSVD and NSTHOSVD decreases their number to 1\% and 2\%, and their norms almost 400 and 1000 times, respectively (we show the convergence curves in Fig.~\ref{fig:Gaussiance_convergence_rate}). Compared to NSTHOSVD, NLRT leaves $1.5$ times more negative elements with a 3 times higher norm. Sketching accelerates NTTSVD and NSTHOSVD 2-3 times, while deterministic NSTHOSVD itself is 3 times faster than NLRT. The increase of the relative error, due to sketching, is within 5\% for most of the methods. We also compare the 4 NLRT components with the auxiliary low-rank tensor \eqref{eq:nlrt_aux} in Table~\ref{tab:Gaussians_NLRT_details}.

\begin{table}[h]\centering
\begin{tabular}{@{}lccccc@{}}
    \toprule
    Method & \thead{Runnung\\time\\(Second)} & \thead{Relative\\error\\(Frobenius)} & \thead{Relative\\error\\(Chebyshev)} & \thead{Negative\\elements\\(Frobenius)} & \thead{Negative\\elements\\(\%)}\\ \hline
    TTSVD                     & 5.9 & $7.4 \cdot 10^{-2}$ & $1.5 \cdot 10^{-1}$ & 5.3 & 41.0    \\ \hline
    NTTSVD, SVD$_r$           & 963 & $8.7 \cdot 10^{-2}$ & $1.8 \cdot 10^{-1}$ & $1.4 \cdot 10^{-2}$ & 1.1  \\  
    NTTSVD, HMT(1, 40)        & 550 & $8.7 \cdot 10^{-2}$ & $1.8 \cdot 10^{-1}$ & $1.4 \cdot 10^{-2}$ & 1.1 \\  
    NTTSVD, HMT(0, 45)        & 452  & $8.7 \cdot 10^{-2}$ & $1.8 \cdot 10^{-1}$ & $1.4 \cdot 10^{-2}$ & 0.89 \\ 
    NTTSVD, Tropp(38, 100)    & 371  & $9.1 \cdot 10^{-2}$ & $1.7 \cdot 10^{-1}$ & $1.4 \cdot 10^{-2}$ & 0.62 \\
    NTTSVD, Tropp(35, 100)    & 325  & $9.4 \cdot 10^{-2}$ & $1.7 \cdot 10^{-1}$ & $1.4 \cdot 10^{-2}$ & 0.49 \\ \hline
    STHOSVD                   & 3.4 & $2.2 \cdot 10^{-2}$ & $7.7 \cdot 10^{-2}$ & 1.8 & 38.0 \\ \hline
    NSTHOSVD, SVD$_r$         & 670 & $2.6 \cdot 10^{-2}$ & $1.0 \cdot 10^{-1}$ & $1.6 \cdot 10^{-3}$ & 1.8\\
    NSTHOSVD, HMT(1, 24)      & 447 & $2.6 \cdot 10^{-2}$ & $1.0 \cdot 10^{-1}$ & $1.6 \cdot 10^{-3}$ & 1.8\\
    NSTHOSVD, HMT(0, 24)      & 346 & $2.6 \cdot 10^{-2}$ & $1.0 \cdot 10^{-1}$ & $1.7 \cdot 10^{-3}$ & 1.3\\
    NSTHOSVD, Tropp(22, 80)   & 230 & $2.7 \cdot 10^{-2}$ & $1.0 \cdot 10^{-1}$ & $1.7 \cdot 10^{-3}$ & 0.99\\
    NSTHOSVD, Tropp(18, 80)   & 205 & $3.5 \cdot 10^{-2}$ & $8.8 \cdot 10^{-2}$ & $1.8 \cdot 10^{-3}$ & 0.75\\ \hline
    NLRT                      & 2096 & $2.6 \cdot 10^{-2}$ & $1.0 \cdot 10^{-1}$ & $5.1 \cdot 10^{-3}$ & 3.1 \\
    \bottomrule
\end{tabular}
\caption{Comparison of NTTSVD, NSTHOSVD, and NLRT for low-rank nonnegative tensor approximation of a $64 \times 64 \times 64 \times 64$ Gaussian mixture with Tucker ranks $(14,14,14,14)$ and TT-ranks $(10,20,10)$: running times, relative errors, and negative elements after 200 iterations.}
\label{tab:Gaussian_timings_errors}
\end{table}

\begin{table}[h]\centering
\begin{tabular}{@{}llcccc@{}}
    \toprule
    Tensor & \thead{Relative\\error\\(Frobenius)} & \thead{Relative\\error\\(Chebyshev)} & \thead{Negative\\elements\\(Frobenius)} & \thead{Negative\\elements\\(Chebyshev)} & \thead{Negative\\elements\\(\%)}\\ \hline
    $\bm{X}_1$ & $2.6 \cdot 10^{-2}$ & $1.0 \cdot 10^{-1}$ & $3.1 \cdot 10^{-3}$ & $2.4 \cdot 10^{-4}$ & 2.0\\
    $\bm{X}_2$ & $2.6 \cdot 10^{-2}$ & $1.0 \cdot 10^{-1}$ & $3.2 \cdot 10^{-3}$ & $2.2 \cdot 10^{-4}$ & 1.9\\
    $\bm{X}_3$ & $2.6 \cdot 10^{-2}$ & $1.0 \cdot 10^{-1}$ & $2.5 \cdot 10^{-3}$ & $1.4 \cdot 10^{-4}$ & 2.7\\
    $\bm{X}_4$ & $2.6 \cdot 10^{-2}$ & $1.0 \cdot 10^{-1}$ & $2.8 \cdot 10^{-3}$ & $2.3 \cdot 10^{-4}$ & 2.2\\
    $\hat{\bm{X}}$ & $2.6 \cdot 10^{-2}$ & $1.0 \cdot 10^{-1}$ & $5.1 \cdot 10^{-3}$ & $3.5 \cdot 10^{-4}$ & 3.1\\
    \bottomrule
\end{tabular}
\caption{Comparison of the NLRT components $\{ \bm{X}_k \}_{k = 1}^{4}$ and the auxiliary tensor $\hat{\bm{X}}$ in Tucker format for low-rank nonnegative tensor approximation of a $64 \times 64 \times 64 \times 64$ Gaussian mixture with Tucker ranks $(14,14,14,14)$: relative errors and negative elements after 200 iterations.}
\label{tab:Gaussians_NLRT_details}
\end{table}

\clearpage
\subsection{Hyperspectral image}\label{subsec:hyperspectral}
Our final example is an openly available hyperspectral image of the Washington DC National Mall\footnote{Data available at \href{https://github.com/JakobSig/HSI2RGB/blob/master/washington_hsi.mat}{\text{https://github.com/JakobSig/HSI2RGB/blob/master/washington\_hsi.mat}}} of size $307 \times 307 \times 191$, the last dimension being the spectral bands (we also linearly scale the elements to $[0,1]$). With Tucker ranks $(40,40,33)$ and TT-ranks $(33,33)$, the image can be compressed 215 and 51 times, respectively. To measure the quality of low-rank approximations, we use 3 values: the relative error in the Frobenius norm, the band-wise mean of the structural similarity index measure (SSIM, \cite{yuan2012hyperspectral}), and the statistical $R^2$ coefficient of determination,
\begin{equation*}
    R^2 = 1 - \frac{\|\bm{X} - \bm{Y}\|_F^2}{\|\bm{X} - \alpha \|_F^2}, \quad \alpha = \frac{1}{\prod\nolimits_{k = 1}^{d} n_k} \sum_{(i_1, \ldots, i_d)} \bm{X}(i_1, \ldots, i_d) \in \mathbb{R},
\end{equation*}
where $\bm{Y}$ is a low-rank approximant.

\begin{figure}[ht]
\begin{subfigure}[b]{1.0\linewidth}
\centering
	\includegraphics[width=\linewidth]{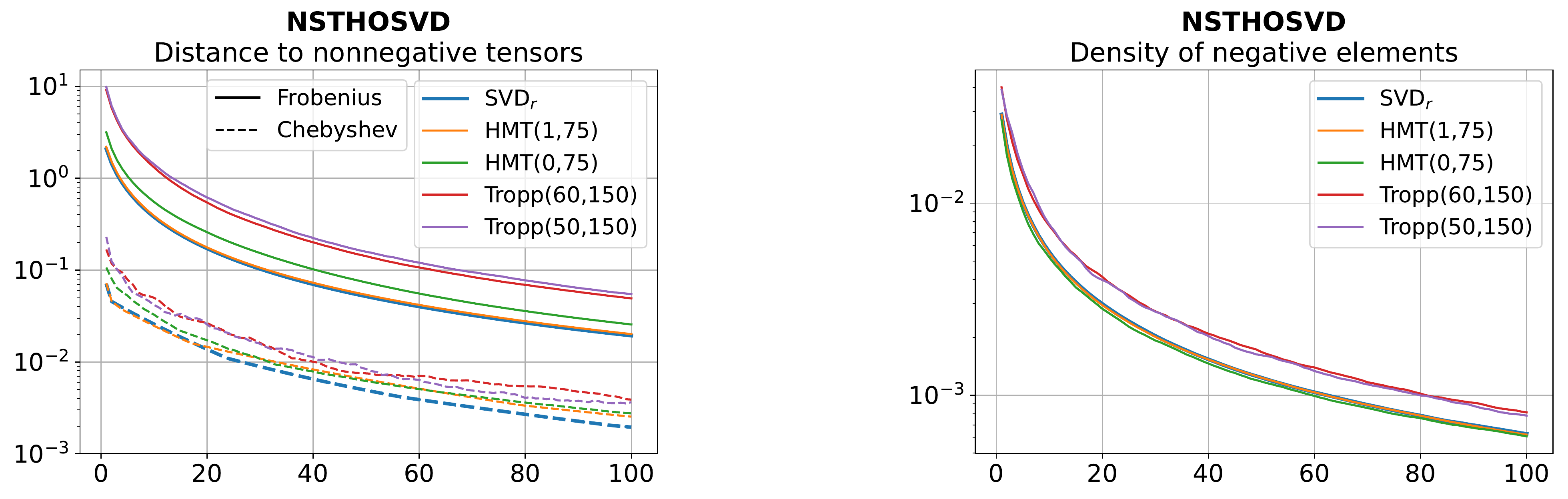}
	\caption{}
\end{subfigure}
\begin{subfigure}[b]{1.0\linewidth}
\centering
	\includegraphics[width=\linewidth]{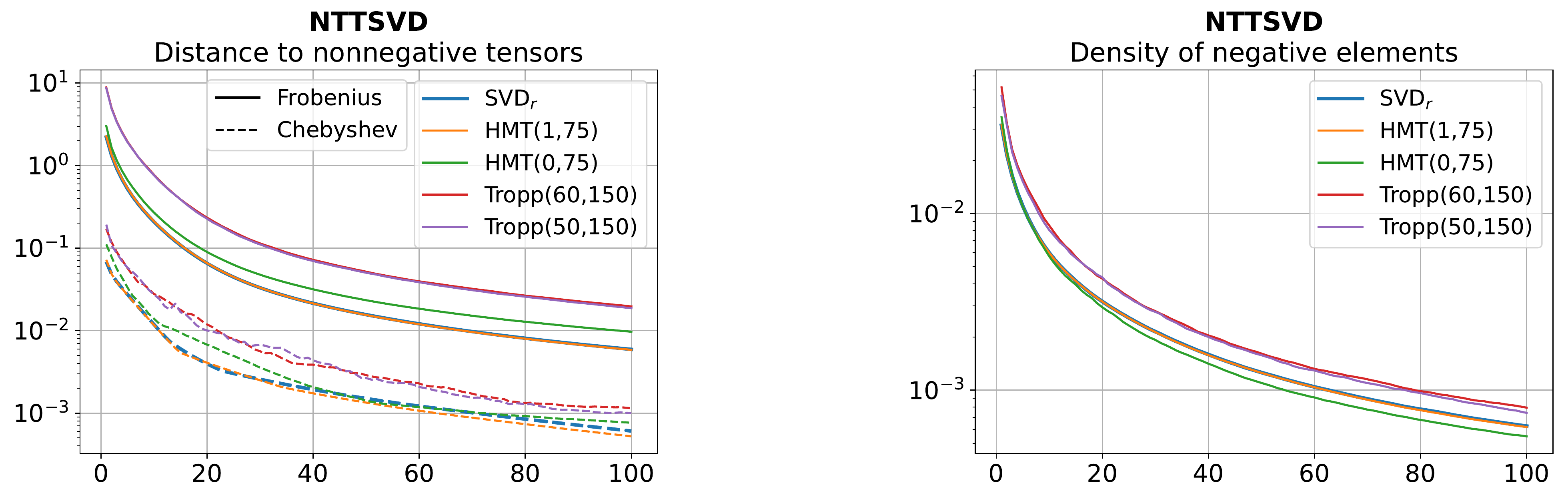}
	\caption{}
\end{subfigure}
\caption{Comparison of deterministic and randomized variants of NSTHOSVD~(a) and NTTSVD~(b) for low-rank nonnegative tensor approximation of a $307 \times 307 \times 191$ hyperspectral image of the Washington DC National Mall with Tucker ranks $(40,40,33)$ and TT-ranks $(33,33)$: the Frobenius and Chebyshev norms of the negative part and the density of negative elements over 100 iterations.}
\label{fig:Convergence_Washington}
\end{figure}

\begin{figure}[ht]
\centering
\includegraphics[width=0.9\linewidth]{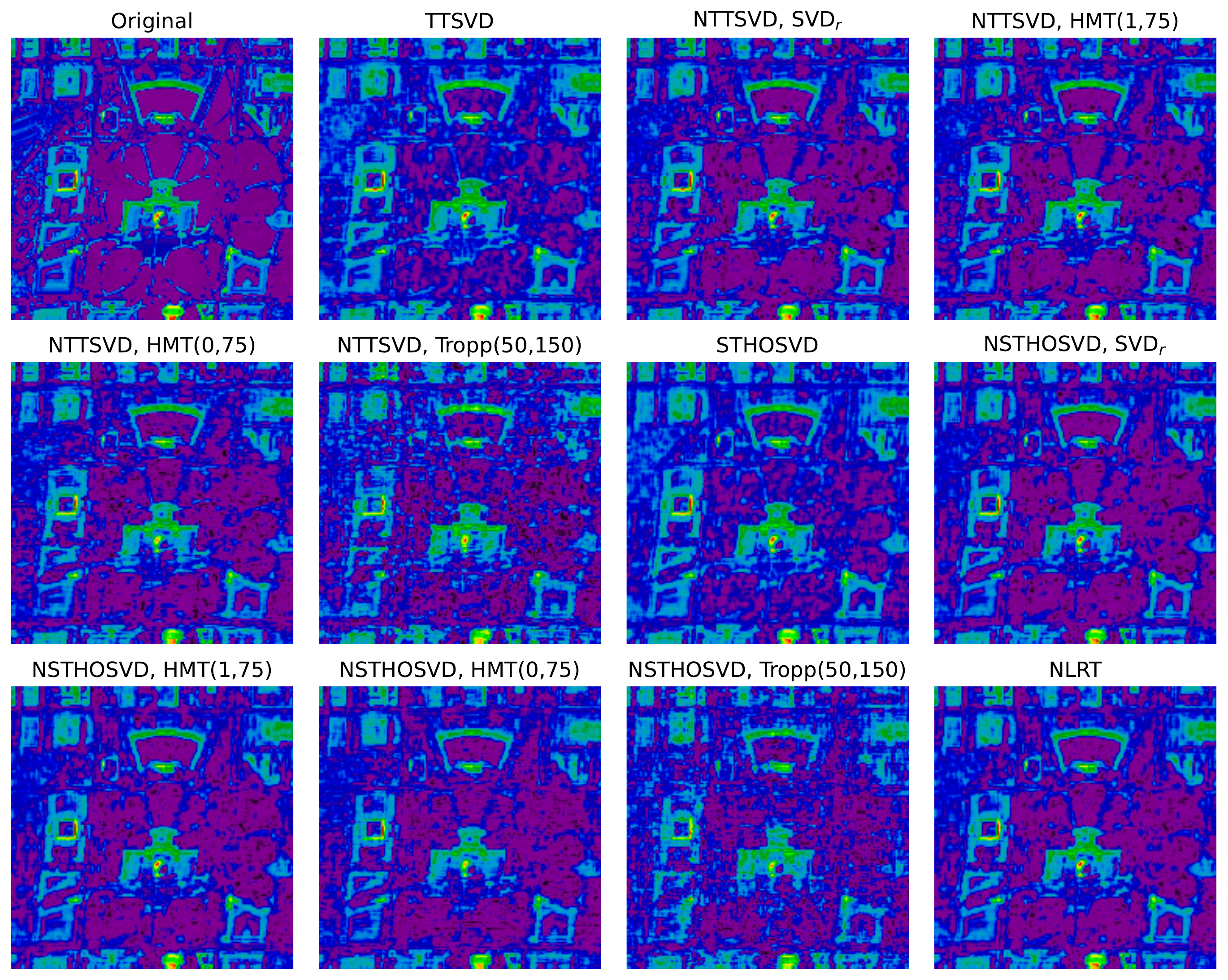}
\caption{Comparison of the approximations of a $307 \times 307 \times 191$ hyperspectral image of the Washington DC National Mall achieved with TTSVD, STHOSVD, and different iterative LRNTA approaches (after 100 iterations). We present the 50th spectral band.}
\label{fig:Visual_Images_Compression_Losses}
\end{figure}

In Table~\ref{tab:Washington_table}, we see that 100 iterations of deterministic NTTSVD and NSTHOSVD lower the Frobenius norm of the negative elements 350 and 100 times, respectively, compared to TTSVD and STHOSVD. Both the relative error and the $R^2$ score stay the same, and SSIM undergoes a reduction by $0.03-0.04$. The randomized variants based on HMT with 1 power-method iteration lead to similar results, but achieve them about 1.5 times faster. With the more computationally efficient randomized approaches, SSIM seems to degrade more severely than the 2 other quality measures. Find the convergence curves for NTTSVD and NSTHOSVD in Fig.~\ref{fig:Convergence_Washington}. The approximation quality achieved with NLRT is identical to deterministic NSTHOSVD; however, the norm of the negative elements is 2.4 times higher, and it runs 4 times slower (6 times compared to randomized NSTHOSVD). As Table~\ref{tab:Washington_NRLT_details} shows, the negative elements in the low-rank auxiliary tensor \eqref{eq:nlrt_aux} have a larger norm than the unfoldings, which NLRT operates on. Finally, in Fig.~\ref{fig:Visual_Images_Compression_Losses} we present the actual images for visual evaluation.

\begin{table}[h]\centering
\begin{tabular}{@{}lccccc@{}}
    \toprule
    Method & \thead{Running\\time\\(Second)} & \thead{Relative\\error\\(Frobenius)} & \thead{SSIM} & $R^2$ & \thead{Negative\\elements\\(Frobenius)}\\ \hline
    TTSVD                    & 3.2 & $1.8 \cdot 10^{-1}$ & 0.66 & 0.94 & $2.2$ \\ \hline
    NTTSVD, SVD$_r$          & 341 & $1.8 \cdot 10^{-1}$ & 0.63 & 0.94 & $6.0 \cdot 10^{-3}$ \\
    NTTSVD, HMT(1, 75)       & 236 & $1.8 \cdot 10^{-1}$ & 0.62 & 0.94 & $6.0 \cdot 10^{-3}$ \\
    NTTSVD, HMT(0, 75)       & 173 & $2.1 \cdot 10^{-1}$ & 0.57 & 0.92 & $1.0 \cdot 10^{-2}$ \\
    NTTSVD, Tropp(60, 150)   & 151 & $2.6 \cdot 10^{-1}$ & 0.46 & 0.88 & $1.8 \cdot 10^{-2}$ \\
    NTTSVD, Tropp(50, 150)   & 102 & $2.7 \cdot 10^{-1}$ & 0.45 & 0.87 & $1.9 \cdot 10^{-2}$ \\ \hline
    STHOSVD                  & 3.4 & $1.8 \cdot 10^{-1}$ & 0.64 & 0.94 & $2.1$ \\ \hline
    NSTHOSVD, SVD$_r$        & 440 & $1.8 \cdot 10^{-1}$ & 0.60 & 0.94 & $1.9 \cdot 10^{-2}$ \\
    NSTHOSVD, HMT(1, 75)     & 296 & $1.9 \cdot 10^{-1}$ & 0.60 & 0.94 & $1.9 \cdot 10^{-2}$ \\
    NSTHOSVD, HMT(0, 75)     & 268 & $2.1 \cdot 10^{-1}$ & 0.53 & 0.92 & $2.6 \cdot 10^{-2}$ \\
    NSTHOSVD, Tropp(60, 150) & 225 & $2.8 \cdot 10^{-1}$ & 0.40 & 0.86 & $5.4 \cdot 10^{-2}$ \\
    NSTHOSVD, Tropp(50, 150) & 131 & $2.9 \cdot 10^{-1}$ & 0.39 & 0.85 & $5.1 \cdot 10^{-2}$ \\
    \hline
    NLRT                     & 1874 & $1.8 \cdot 10^{-1}$ & 0.60 & 0.94 & $4.6 \cdot 10^{-2}$ \\
    \bottomrule
\end{tabular}
\caption{Comparison of NTTSVD, NSTHOSVD, and NLRT for low-rank nonnegative tensor approximation of a $307 \times 307 \times 191$ hyperspectral image of the Washington DC National Mall with Tucker ranks $(40,40,33)$ and TT-ranks $(33,33)$: running times, relative errors, negative elements, SSIM, and $R^2$ score after 100 iterations.}
\label{tab:Washington_table}
\end{table}

\begin{table}[h]\centering
\begin{tabular}{@{}lcccc@{}}
    \toprule
    Tensor & \thead{Relative error\\(Frobenius)} & \thead{SSIM} & $R^2$ & \thead{Negative elements\\(Frobenius)} \\ \hline
    $\bm{X}_1$     & $1.8 \cdot 10^{-1}$ & 0.60 & 0.94 & $3.6 \cdot 10^{-2}$ \\
    $\bm{X}_2$     & $1.8 \cdot 10^{-1}$ & 0.60 & 0.94 & $3.6 \cdot 10^{-2}$ \\
    $\bm{X}_3$     & $1.8 \cdot 10^{-1}$ & 0.60 & 0.94 & $1.5 \cdot 10^{-2}$ \\
    $\hat{\bm{X}}$ & $1.8 \cdot 10^{-1}$ & 0.60 & 0.94 & $4.6 \cdot 10^{-2}$ \\
    \bottomrule
\end{tabular}
\caption{Comparison of the NLRT components $\{ \bm{X}_k \}_{k = 1}^{3}$ and the auxiliary tensor $\hat{\bm{X}}$ in Tucker format for low-rank nonnegative tensor approximation of a $307 \times 307 \times 191$ hyperspectral image of the Washington DC National Mall with Tucker ranks $(40,40,33)$: negative elements, relative errors, SSIM, and $R^2$ score after 100 iterations.}
\label{tab:Washington_NRLT_details}
\end{table}

\clearpage
\section{Conclusion}\label{sec:conclusions}
\begin{figure}[t]
\centering
\begin{subfigure}[b]{0.3\textwidth}
	\includegraphics[width=1.0\linewidth]{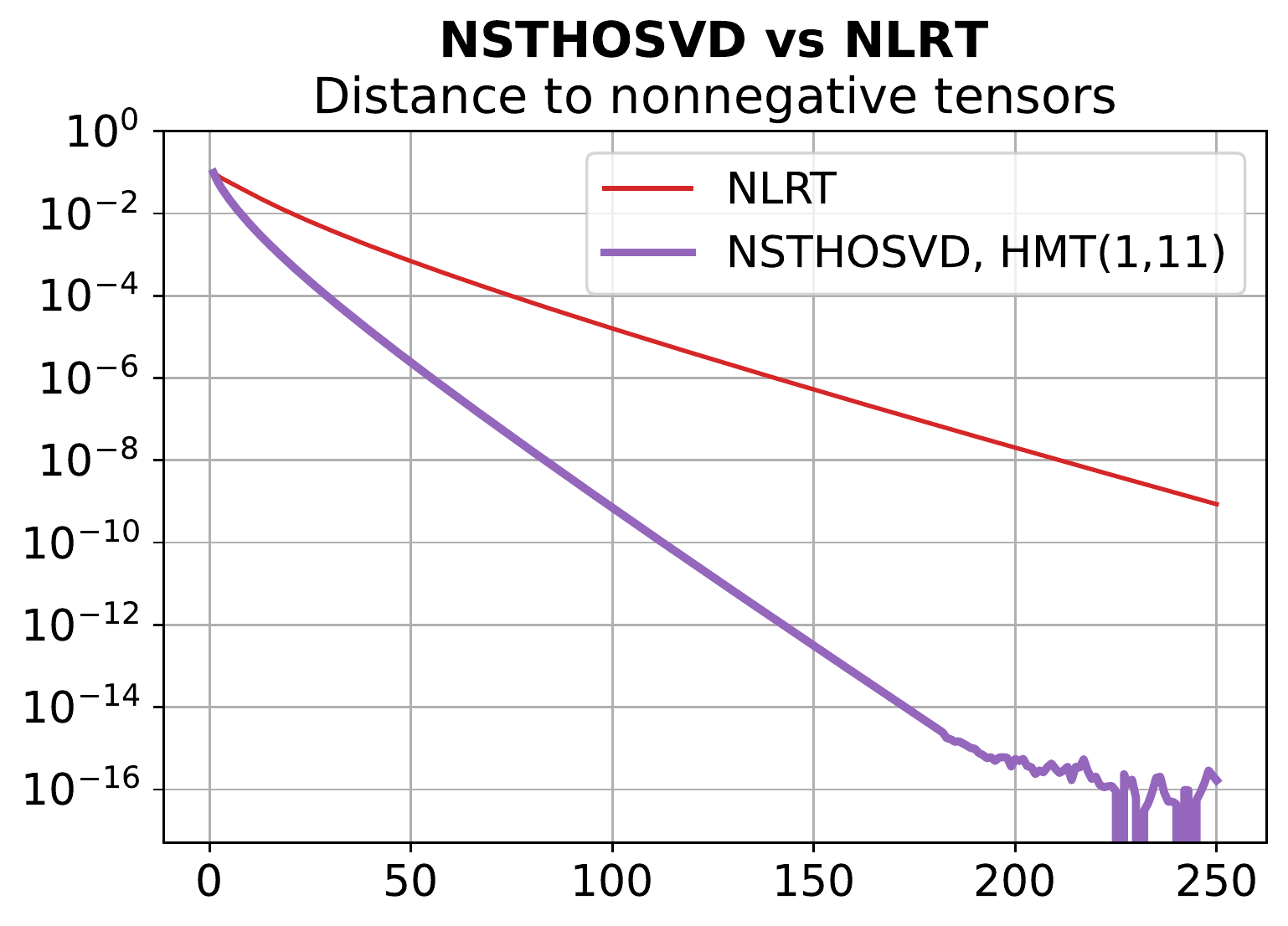}
\caption{}
\end{subfigure}\hfill%
\begin{subfigure}[b]{0.3\textwidth}
	\includegraphics[width=1.0\linewidth]{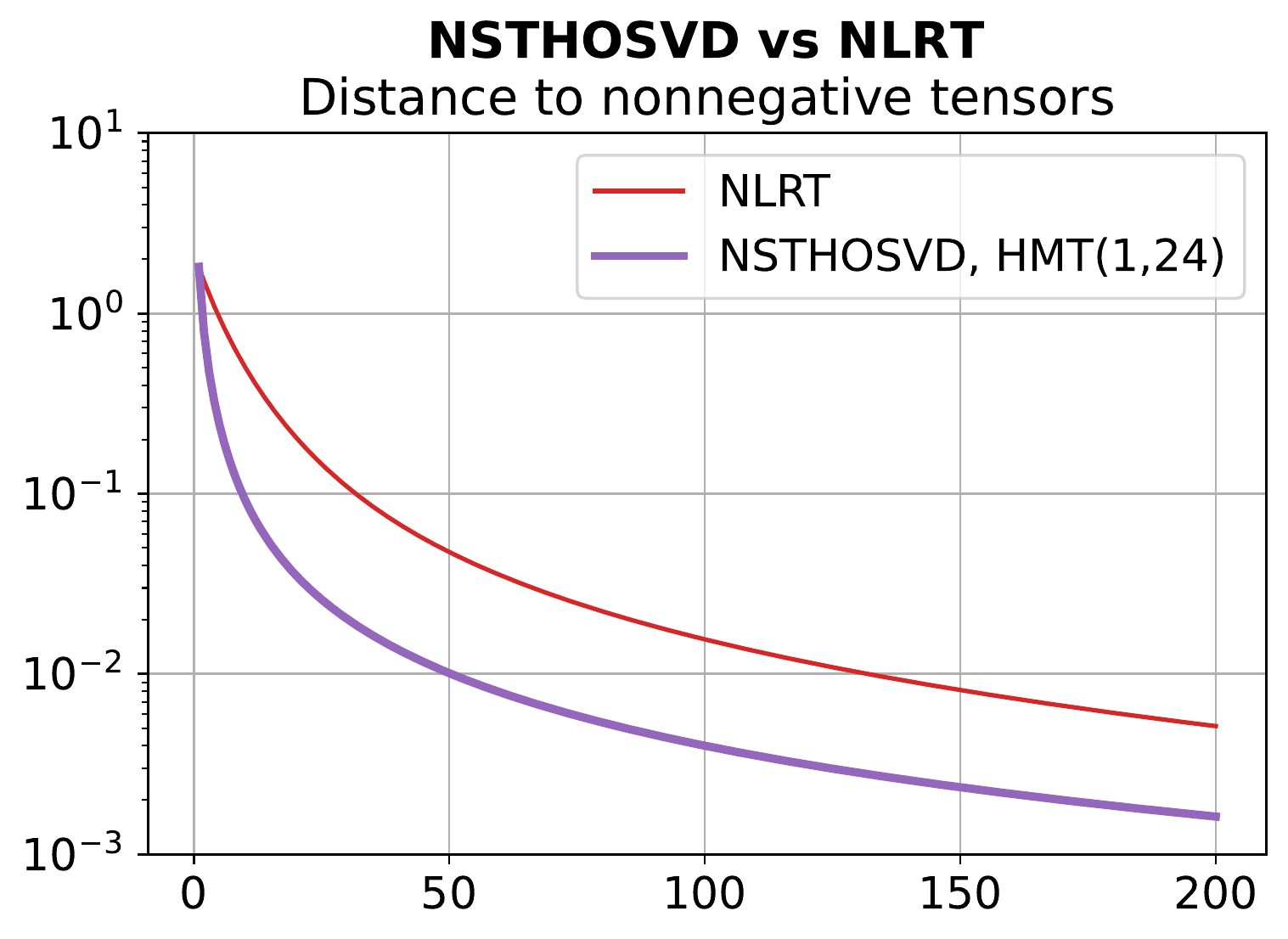}
\caption{}
\end{subfigure}\hfill%
\begin{subfigure}[b]{0.3\textwidth}
	\includegraphics[width=1.0\linewidth]{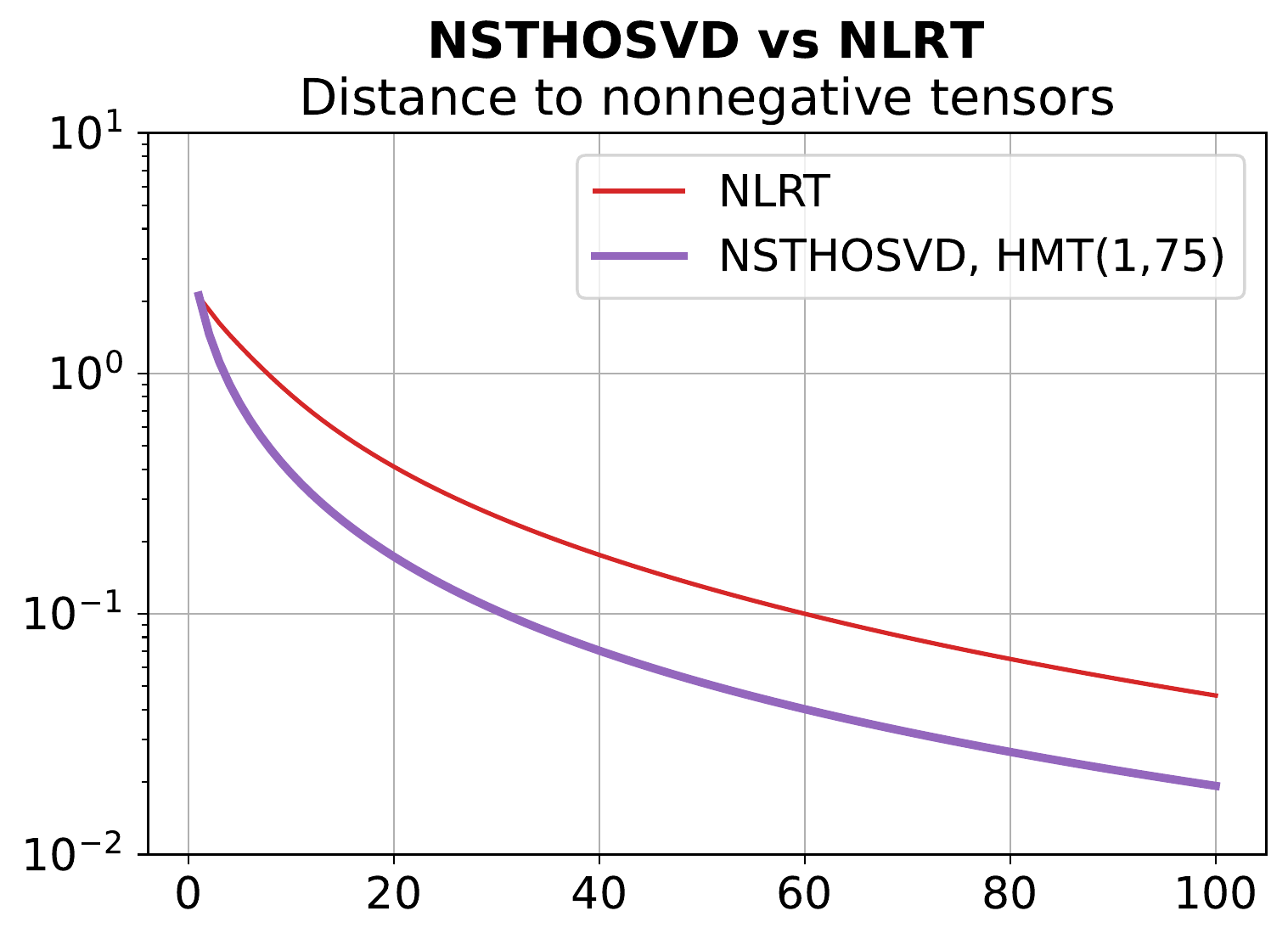}
\caption{}
\end{subfigure}
\caption{Comparison of randomized NSTHOSVD and NLRT for low-rank nonnegative tensor approximation: the decay of the Frobenius norm of the negative elements for the Hilbert tensor~(a), the Gaussian mixture~(b), and the hyperspectral image~(c).}
\label{fig:All_NRLT_distance}
\end{figure}

In this work, we looked at a natural multidimensional extension of randomized alternating projections for the LRNMA problem \cite{matveev2022sketching} and proposed two algorithms, NSTHOSVD and NTTSVD, for the Tucker and tensor train formats, respectively. The numerical experiments showed that our approach allows to reduce the number (and the absolute value) of the negative elements in the low-rank approximation without significant loss of accuracy. Comparing with the NLRT method \cite{jiang2020nonnegative}, which was developed for the Tucker case, we observed that our algorithm NSTHOSVD is superior in terms of computational efficiency per iteration and in how fast it reduces the negative elements (see Fig.~\ref{fig:All_NRLT_distance}).

The use of randomization allowed us to obtain algorithms, whose complexity scales linearly with the number of elements of the tensor, thereby achieving balance in complexity of low-rank and nonnegative projections. Moreover, by choosing the configuration parameters of randomized sketching (such as the oversampling, the distribution of the random matrix) one can tune the methods to achieve the desired trade-off between speed and accuracy.

Though the proposed algorithms work in numerical experiments, they still require a proof of convergence. We will study their theoretical properties in future papers.

\section*{Acknowledgements}
This work was supported by Russian Science Foundation (project 21-71-10072). Sergey Matveev is grateful to Sirius University for opportunity to finalize a significant part of his work on this paper during the summer school ``Matrix methods and modelling in geophysics and life sciences''. All authors have equally contributed to the study conception and design of the research present in this article.

\section*{Data Availability}

The datasets generated during and/or analysed during the current study are available in the github repository \url{https://github.com/azamat11235/NLRTA}.

\bibliography{submitted_version/references}

\begin{thebibliography}{10}

\bibitem{cichocki2017tensor}
A.~Cichocki, A.-H. Phan, Q.~Zhao, N.~Lee, I.~Oseledets, M.~Sugiyama, D.~P.
  Mandic, {\em et~al.}, ``Tensor networks for dimensionality reduction and
  large-scale optimization: Part 2 applications and future perspectives,'' {\em
  Foundations and Trends{\textregistered} in Machine Learning}, vol.~9, no.~6,
  pp.~431--673, 2017.

\bibitem{sidiropoulos2017tensor}
N.~D. Sidiropoulos, L.~De~Lathauwer, X.~Fu, K.~Huang, E.~E. Papalexakis, and
  C.~Faloutsos, ``Tensor decomposition for signal processing and machine
  learning,'' {\em IEEE Transactions on Signal Processing}, vol.~65, no.~13,
  pp.~3551--3582, 2017.

\bibitem{khoromskij2018tensor}
B.~N. Khoromskij, ``Tensor numerical methods in scientific computing,'' in {\em
  Tensor Numerical Methods in Scientific Computing}, De Gruyter, 2018.

\bibitem{kazeev2014direct}
V.~Kazeev, M.~Khammash, M.~Nip, and C.~Schwab, ``Direct solution of the
  chemical master equation using quantized tensor trains,'' {\em PLoS
  computational biology}, vol.~10, no.~3, p.~e1003359, 2014.

\bibitem{matveev2016tensor}
S.~A. Matveev, D.~A. Zheltkov, E.~E. Tyrtyshnikov, and A.~P. Smirnov, ``Tensor
  train versus {M}onte {C}arlo for the multicomponent {S}moluchowski
  coagulation equation,'' {\em Journal of Computational Physics}, vol.~316,
  pp.~164--179, 2016.

\bibitem{allmannrahn2022parallel}
F.~Allmann-Rahn, R.~Grauer, and K.~Kormann, ``A parallel low-rank solver for
  the six-dimensional vlasov-maxwell equations,'' {\em arXiv preprint
  arXiv:2201.03471}, 2022.

\bibitem{liu2012tensor}
J.~Liu, P.~Musialski, P.~Wonka, and J.~Ye, ``Tensor completion for estimating
  missing values in visual data,'' {\em IEEE transactions on pattern analysis
  and machine intelligence}, vol.~35, no.~1, pp.~208--220, 2012.

\bibitem{budzinskiy2021tensor}
S.~Budzinskiy and N.~Zamarashkin, ``Tensor train completion: local recovery
  guarantees via {R}iemannian optimization,'' {\em arXiv preprint
  arXiv:2110.03975}, 2021.

\bibitem{zheltkov2020global}
D.~Zheltkov and E.~Tyrtyshnikov, ``Global optimization based on
  {T}{T}-decomposition,'' {\em Russian Journal of Numerical Analysis and
  Mathematical Modelling}, vol.~35, no.~4, pp.~247--261, 2020.

\bibitem{yang2017tensor}
Y.~Yang, D.~Krompass, and V.~Tresp, ``Tensor-train recurrent neural networks
  for video classification,'' in {\em International Conference on Machine
  Learning}, pp.~3891--3900, PMLR, 2017.

\bibitem{dolgov2020approximation}
S.~Dolgov, K.~Anaya-Izquierdo, C.~Fox, and R.~Scheichl, ``Approximation and
  sampling of multivariate probability distributions in the tensor train
  decomposition,'' {\em Statistics and Computing}, vol.~30, no.~3,
  pp.~603--625, 2020.

\bibitem{novikov2021tensor}
G.~S. Novikov, M.~E. Panov, and I.~V. Oseledets, ``Tensor-train density
  estimation,'' in {\em Uncertainty in Artificial Intelligence},
  pp.~1321--1331, PMLR, 2021.

\bibitem{hur2022generative}
Y.~Hur, J.~G. Hoskins, M.~Lindsey, E.~M. Stoudenmire, and Y.~Khoo, ``Generative
  modeling via tensor train sketching,'' {\em arXiv preprint arXiv:2202.11788},
  2022.

\bibitem{shcherbakova2019nonnegative1}
E.~Shcherbakova and E.~Tyrtyshnikov, ``Nonnegative tensor train factorizations
  and some applications,'' in {\em International Conference on Large-Scale
  Scientific Computing}, pp.~156--164, Springer, 2019.

\bibitem{manzini2021nonnegative}
G.~Manzini, E.~Skau, D.~P. Truong, and R.~Vangara, ``Nonnegative tensor-train
  low-rank approximations of the {S}moluchowski coagulation equation,'' in {\em
  International Conference on Large-Scale Scientific Computing}, pp.~342--350,
  Springer, 2021.

\bibitem{cichocki2009nonnegative}
A.~Cichocki, R.~Zdunek, A.~H. Phan, and S.-i. Amari, {\em Nonnegative matrix
  and tensor factorizations: applications to exploratory multi-way data
  analysis and blind source separation}.
\newblock John Wiley \& Sons, 2009.

\bibitem{leplat2020blind}
V.~Leplat, N.~Gillis, and A.~M. Ang, ``Blind audio source separation with
  minimum-volume beta-divergence {N}{M}{F},'' {\em IEEE Transactions on Signal
  Processing}, vol.~68, pp.~3400--3410, 2020.

\bibitem{gillis2020nonnegative}
N.~Gillis, {\em Nonnegative matrix factorization}.
\newblock SIAM, 2020.

\bibitem{vanluyten2008nonnegative}
B.~Vanluyten, J.~C. Willems, and B.~De~Moor, ``Nonnegative matrix factorization
  without nonnegativity constraints on the factors,'' {\em Submitted for
  publication}, 2008.

\bibitem{grussler2015optimal}
C.~Grussler and A.~Rantzer, ``On optimal low-rank approximation of non-negative
  matrices,'' in {\em 2015 54th IEEE Conference on Decision and Control (CDC)},
  pp.~5278--5283, IEEE, 2015.

\bibitem{grussler2018low}
C.~Grussler, A.~Rantzer, and P.~Giselsson, ``Low-rank optimization with convex
  constraints,'' {\em IEEE Transactions on Automatic Control}, vol.~63, no.~11,
  pp.~4000--4007, 2018.

\bibitem{andersson2017convex}
F.~Andersson, M.~Carlsson, and C.~Olsson, ``Convex envelopes for fixed rank
  approximation,'' {\em Optimization Letters}, vol.~11, no.~8, pp.~1783--1795,
  2017.

\bibitem{song2020nonnegative}
G.-J. Song and M.~K. Ng, ``Nonnegative low rank matrix approximation for
  nonnegative matrices,'' {\em Applied Mathematics Letters}, vol.~105,
  p.~106300, 2020.

\bibitem{song2020tangent}
G.~Song, M.~K. Ng, and T.-X. Jiang, ``Tangent space based alternating
  projections for nonnegative low rank matrix approximation,'' {\em arXiv
  preprint arXiv:2009.03998}, 2020.

\bibitem{matveev2022sketching}
S.~A. Matveev and S.~Budzinskiy, ``Sketching for low-rank nonnegative matrix
  approximation: a numerical study,'' {\em arXiv preprint arXiv:2201.11154},
  2022.

\bibitem{zhu2021approximate}
H.~Zhu, M.~K. Ng, and G.-J. Song, ``An approximate augmented {L}agrangian
  method for nonnegative low-rank matrix approximation,'' {\em Journal of
  Scientific Computing}, vol.~88, no.~2, pp.~1--22, 2021.

\bibitem{tucker1964extension}
L.~R. Tucker {\em et~al.}, ``The extension of factor analysis to
  three-dimensional matrices,'' {\em Contributions to mathematical psychology},
  vol.~110119, 1964.

\bibitem{oseledets2011tensor}
I.~V. Oseledets, ``Tensor-train decomposition,'' {\em SIAM Journal on
  Scientific Computing}, vol.~33, no.~5, pp.~2295--2317, 2011.

\bibitem{jiang2020nonnegative}
T.-X. Jiang, M.~K. Ng, J.~Pan, and G.~Song, ``Nonnegative low rank tensor
  approximation and its application to multi-dimensional images,'' {\em arXiv
  preprint arXiv:2007.14137}, 2020.

\bibitem{kolda2009tensor}
T.~G. Kolda and B.~W. Bader, ``Tensor decompositions and applications,'' {\em
  SIAM review}, vol.~51, no.~3, pp.~455--500, 2009.

\bibitem{de2000multilinear}
L.~De~Lathauwer, B.~De~Moor, and J.~Vandewalle, ``A multilinear singular value
  decomposition,'' {\em SIAM journal on Matrix Analysis and Applications},
  vol.~21, no.~4, pp.~1253--1278, 2000.

\bibitem{dongarra2018singular}
J.~Dongarra, M.~Gates, A.~Haidar, J.~Kurzak, P.~Luszczek, S.~Tomov, and
  I.~Yamazaki, ``The singular value decomposition: Anatomy of optimizing an
  algorithm for extreme scale,'' {\em SIAM review}, vol.~60, no.~4,
  pp.~808--865, 2018.

\bibitem{oseledets2009breaking}
I.~V. Oseledets and E.~E. Tyrtyshnikov, ``Breaking the curse of dimensionality,
  or how to use {S}{V}{D} in many dimensions,'' {\em SIAM Journal on Scientific
  Computing}, vol.~31, no.~5, pp.~3744--3759, 2009.

\bibitem{martinsson2020randomized}
P.-G. Martinsson and J.~A. Tropp, ``Randomized numerical linear algebra:
  Foundations and algorithms,'' {\em Acta Numerica}, vol.~29, pp.~403--572,
  2020.

\bibitem{halko2011finding}
N.~Halko, P.-G. Martinsson, and J.~A. Tropp, ``Finding structure with
  randomness: Probabilistic algorithms for constructing approximate matrix
  decompositions,'' {\em SIAM review}, vol.~53, no.~2, pp.~217--288, 2011.

\bibitem{tropp2017practical}
J.~A. Tropp, A.~Yurtsever, M.~Udell, and V.~Cevher, ``Practical sketching
  algorithms for low-rank matrix approximation,'' {\em SIAM Journal on Matrix
  Analysis and Applications}, vol.~38, no.~4, pp.~1454--1485, 2017.

\bibitem{GolubVanLoanMatrix2013}
G.~H. Golub and C.~F. Van~Loan, {\em Matrix Computations}.
\newblock Johns {{Hopkins}} Studies in the Mathematical Sciences, {Baltimore}:
  {The Johns Hopkins University Press}, fourth edition~ed., 2013.

\bibitem{vannieuwenhoven2012new}
N.~Vannieuwenhoven, R.~Vandebril, and K.~Meerbergen, ``A new truncation
  strategy for the higher-order singular value decomposition,'' {\em SIAM
  Journal on Scientific Computing}, vol.~34, no.~2, pp.~A1027--A1052, 2012.

\bibitem{che2019randomized}
M.~Che and Y.~Wei, ``Randomized algorithms for the approximations of {T}ucker
  and the tensor train decompositions,'' {\em Advances in Computational
  Mathematics}, vol.~45, no.~1, pp.~395--428, 2019.

\bibitem{ahmadi2021randomized}
S.~Ahmadi-Asl, S.~Abukhovich, M.~G. Asante-Mensah, A.~Cichocki, A.~H. Phan,
  T.~Tanaka, and I.~Oseledets, ``Randomized algorithms for computation of
  {T}ucker decomposition and higher order {S}{V}{D} ({H}{O}{S}{V}{D}),'' {\em
  IEEE Access}, vol.~9, pp.~28684--28706, 2021.

\bibitem{huber2017randomized}
B.~Huber, R.~Schneider, and S.~Wolf, ``A randomized tensor train singular value
  decomposition,'' in {\em Compressed sensing and its applications},
  pp.~261--290, Springer, 2017.

\bibitem{kressner2022streaming}
D.~Kressner, B.~Vandereycken, and R.~Voorhaar, ``Streaming tensor train
  approximation,'' {\em arXiv preprint arXiv:2208.02600}, 2022.

\bibitem{escalante2011alternating}
R.~Escalante and M.~Raydan, {\em Alternating projection methods}.
\newblock SIAM, 2011.

\bibitem{deutsch1984rate}
F.~Deutsch, ``Rate of convergence of the method of alternating projections,''
  in {\em Parametric optimization and approximation}, pp.~96--107, Springer,
  1984.

\bibitem{kayalar1988error}
S.~Kayalar and H.~L. Weinert, ``Error bounds for the method of alternating
  projections,'' {\em Mathematics of Control, Signals and Systems}, vol.~1,
  no.~1, pp.~43--59, 1988.

\bibitem{bauschke1994dykstra}
H.~H. Bauschke and J.~M. Borwein, ``Dykstra's alternating projection algorithm
  for two sets,'' {\em Journal of Approximation Theory}, vol.~79, no.~3,
  pp.~418--443, 1994.

\bibitem{bauschke1993convergence}
H.~H. Bauschke and J.~M. Borwein, ``On the convergence of von neumann's
  alternating projection algorithm for two sets,'' {\em Set-Valued Analysis},
  vol.~1, no.~2, pp.~185--212, 1993.

\bibitem{deutsch2001best}
F.~Deutsch, {\em Best approximation in inner product spaces}, vol.~7.
\newblock Springer, 2001.

\bibitem{poliquin2000local}
R.~Poliquin, R.~Rockafellar, and L.~Thibault, ``Local differentiability of
  distance functions,'' {\em Transactions of the American mathematical
  Society}, vol.~352, no.~11, pp.~5231--5249, 2000.

\bibitem{lewis2008alternating}
A.~S. Lewis and J.~Malick, ``Alternating projections on manifolds,'' {\em
  Mathematics of Operations Research}, vol.~33, no.~1, pp.~216--234, 2008.

\bibitem{lewis2009local}
A.~S. Lewis, D.~R. Luke, and J.~Malick, ``Local linear convergence for
  alternating and averaged nonconvex projections,'' {\em Foundations of
  Computational Mathematics}, vol.~9, no.~4, pp.~485--513, 2009.

\bibitem{andersson2013alternating}
F.~Andersson and M.~Carlsson, ``Alternating projections on nontangential
  manifolds,'' {\em Constructive approximation}, vol.~38, no.~3, pp.~489--525,
  2013.

\bibitem{uschmajew2020geometric}
A.~Uschmajew and B.~Vandereycken, ``Geometric methods on low-rank matrix and
  tensor manifolds,'' in {\em Handbook of variational methods for nonlinear
  geometric data}, pp.~261--313, Springer, 2020.

\bibitem{luke2013prox}
D.~R. Luke, ``Prox-regularity of rank constraint sets and implications for
  algorithms,'' {\em Journal of Mathematical Imaging and Vision}, vol.~47,
  no.~3, pp.~231--238, 2013.

\bibitem{yuan2012hyperspectral}
Q.~Yuan, L.~Zhang, and H.~Shen, ``Hyperspectral image denoising employing a
  spectral--spatial adaptive total variation model,'' {\em IEEE Transactions on
  Geoscience and Remote Sensing}, vol.~50, no.~10, pp.~3660--3677, 2012.

\end{thebibliography}
\bibliographystyle{ieeetr}
\addcontentsline{toc}{section}{References}

\end{document}